\documentclass{article}

\usepackage{PRIMEarxiv}

\usepackage[utf8]{inputenc} % allow utf-8 input
\usepackage[T1]{fontenc}    % use 8-bit T1 fonts
\usepackage{hyperref}       % hyperlinks
\usepackage{url}  
\usepackage{booktabs}       % professional-quality tables
\usepackage{amsfonts}       % blackboard math symbols
\usepackage{nicefrac}       % compact symbols for 1/2, etc.
\usepackage{microtype}      % microtypography
\usepackage{lipsum}
\usepackage{fancyhdr}       % header
\usepackage{graphicx}       % graphics
\usepackage{amsmath}
\usepackage{cleveref}

\newenvironment{proof}[1][Proof]{\par\noindent\text{#1. }\ }{\hfill $\square$\par}

\graphicspath{{media/}}     % organize your images and other figures under media/ folder

\usepackage{color}
\usepackage{booktabs}
\usepackage{subcaption}
\usepackage{float}

\newtheorem{theorem}{Theorem}
\newtheorem{lemma}{Lemma}
\newtheorem{definition}{Definition}
\newtheorem{proposition}{Proposition}
\newtheorem{corollary}{Corollary}
\newtheorem{remark}{Remark}

\title{Optimizing Impulsive Releases: A Species Competition Model
}

\author{
  Jéssica C.S. Alves\textsuperscript{1}, Sergio M. Oliva \\
  University of São Paulo\\ Institute of Mathematics and Statistics\\ Department of Applied Mathematics\\ São Paulo, SP, 05508-090, Brazil.\\
  \texttt{alvesj@ime.usp.br, soliva@usp.br} \\
   \And
  Christian E. Schaerer \\
  National University of Asunción\\ Polytechnic School\\ Campus UNA\\ San Lorenzo, Central, P.O. Box 2111 SL, Paraguay.\\
  \texttt{cschaer@pol.una.py} \\
}

\begin{document}
\maketitle

\begin{abstract}
This study focuses on optimizing species release $S_2$ to control species population $S_1$ through impulsive release strategies. We investigate the conditions required to remove species $S_1$, which is equivalent to the establishment of $S_2$. The research includes a theoretical analysis that examines the positivity, existence, and uniqueness of solutions, the conditions ensuring global stability, and a sufficient condition for controlling the $S_1$-free solution. In addition, we formulate an optimal control problem to maximize the effectiveness of $S_2$ releases, manage the population of $S_1$, and minimize the costs associated with this intervention strategy. Numerical simulations are conducted to validate the proposed theories and allow visualization of population dynamics under various release scenarios.
\end{abstract}

% keywords can be removed

\section{Introduction}
\footnotetext[1]{Corresponding author: alvesj@ime.usp.br.}
Population control in ecological systems is a highly relevant issue, especially when two or more species (or phenotypes of the same species) directly compete for limited resources. These interactions can lead to competitive exclusion, in which only one species survives in the long term. In many cases, managing such interactions involves introducing a new species into a region already occupied by another, aiming to control or suppress the original population \cite{HOLT2017Species}.

Various population control methods are used, ranging from chemical \cite{brunner1994integratedChemical1,carson2023overcomingI3,gautam2023pesticideChemical2,lees2023insecticidesI1} and mechanical \cite{adhikari2022insectMechanical1,vincent2009physicalM2} strategies to the use of biological agents \cite{brunner1994integratedChemical1,barbosa2018modellingI2, campo2018optimalW9,onen2023mosquito}. The introduction of a competing species is a widely used approach, as it allows natural ecological interactions, such as competition, to sustainably control the target population \cite{lopes2023exploringW4,almeida2022optimalSterile1,barclay1980sterile7,campo2017optimalW10}.

Competitive interactions between species have been extensively studied in the literature through different models and techniques to represent the dynamics of populations competing for limited resources. Some studies explore competition models that incorporate the Allee effect, highlighting the importance of nonlinear mechanisms in population interactions \cite{Munkaila, Sophia}. Other works use impulsive models to describe dynamics with discontinuous events, such as periodic interventions \cite{YueMeng, Shuwen, Manuel}. Significant contributions have also been made in the context of optimal control, where intervention strategies are formulated under cost and efficacy constraints \cite{Munkaila}. Additionally, feedback-based impulsive control approaches have been proposed to reinforce the applicability of impulsive systems in population management \cite{XuJing}.

In parallel, optimal control techniques have emerged as powerful tools to determine intervention strategies that minimize costs and maximize effectiveness, particularly in ecological and epidemiological contexts \cite{raja2024improved,almeida2019optimal,pei2018optimizingPests1,berkovitz2013optimal,cesari1983optimization}. Moreover, impulsive systems have gained prominence as a natural approach to modeling time-discrete interventions, such as periodic releases of individuals or treatments applied at specific moments, by integrating these three elements: species competition, optimal control, and impulsive systems. The present work aligns with current research trends and contributes to advancing the field by proposing an original formulation aimed optimize intervention strategies in biological systems with discontinuous dynamics.

Many studies use mathematical models to represent interactions between competing or predatory species, often based on continuous approaches \cite{campo2018optimalW9,ogunlade2020modelingC1,silva2020modelingC2,perez2020classC3}. These formulations assume that control or introduction of new species occurs continuously, which, although theoretically feasible, does not always reflect practical reality. In practice, interventions such as the introduction of a competing species generally occur periodically or impulsively due to logistical, financial, and operational constraints. To address these limitations, some studies propose alternative approaches, such as impulsive releases. Examples include the release of mosquitoes with Wolbachia to control the wild Aedes aegypti population \cite{almeida2022vectorImpulsive1,li2024modelingImpulsive3,liu2023analysisImpulsive2,dianavinnarasi2021controllingImpulsive4}, periodic and impulsive release of sterile mosquitoes \cite{huang2021studySterile4,huang2017modellingSterile7,huang2021impulsiveSterile5,li2020impulsiveSterile6}, and similar strategies applied to other insect populations \cite{liu2023analysisImpulsive2,pei2018optimizingPests1,wang2011analysisPests2}.

In this work, we adapt the interaction model between wild and Wolbachia-infected female mosquitoes of the species Aedes aegypti, originally proposed in \cite{Campo2017}, to represent the competition between two species, denoted by $S_1$ and $S_2$. The adapted model incorporates impulsive releases of species $S_2$ and is formulated as a system of impulsive differential equations with jumps occurring at the initial condition of each intervention cycle. This modeling choice provides a realistic representation of biological control practices, in which the introduction of individuals from a competing population typically occurs at discrete time intervals rather than continuously. By placing the impulse at the initial state of each period, the model captures the cumulative effect of releases in a mathematically concise and biologically interpretable way.

One of the main contributions of this study lies in the analytical derivation of a sufficient condition that guarantees both the elimination of species $S_1$ and the global stability of species $S_2$ under an impulsive release strategy. This result not only characterizes the long-term behavior of the system but also provides a clear and actionable threshold for successful intervention. In addition, we formulate and solve an optimal control problem that determines the most efficient release strategy minimizing the number and magnitude of interventions while ensuring the fixation of $S_2$.This combination of impulsive modeling, theoretical contributions, and cost-effective optimization constitutes the novelty of the proposed approach.

The structure of the paper is as follows. In Section 2, we present the formulation of the adapted model from \cite{Campo2017}, including the impulsive differential equations that represent the periodic release of individuals from population $S_2$. We then describe the model components in detail, discussing the parameters involved and the specific conditions required to ensure biologically consistent population dynamics. The section concludes with a theorem on the equilibrium points and their stability, which was proposed in \cite{Campo2017} and is being adapted to the context in this study.

Section 3 is dedicated to the analysis of the model dynamics. We begin by establishing fundamental results concerning the existence, uniqueness, positivity, and boundedness of solutions to the impulsive system, ensuring temporal consistency and well-posedness. Then, we investigate the existence of a solution in which species $S_1$ is eliminated and provide a detailed stability analysis of this solution. In particular, we derive a sufficient condition that guarantees the global stability of the $S_1$-free solution, along with a criterion for the amount of individuals of species $S_2$ required to eliminate species $S_1$ through impulsive releases.

In Section 4, we formulate an optimal control problem aimed at minimizing the total number of individuals released during the intervention interval $[0,T]$, where $T$ denotes the final observation time. At the same time, the strategy must ensure that the population of species $S_1$ is reduced below its survival threshold associated with the Allee effect by time $T$. We also prove the existence of at least one optimal solution to this control problem.

Section 5 presents numerical simulations based on the interaction between two subspecies: wild females and Wolbachia-infected females of Aedes aegypti. First, we simulate the impulsive model to validate the theoretical findings from Section 3, confirming the applicability of the sufficient condition for global stability of the $S_1$-free solution. We then simulate the optimal control problem, using predefined impulse times and varying the final intervention time $T$. These simulations allow us to assess how different parameters affect the effectiveness of control and guide the optimization of release strategies to reduce costs while ensuring the elimination of species $S_1$.

Finally, Section 6 presents the main conclusions of the study. We discuss the relevance of impulsive release strategies in the population dynamics of species $S_1$ and $S_2$, emphasizing the role of optimal control in enhancing the effectiveness and efficiency of ecological interventions. We also propose directions for future research, highlighting opportunities to extend and deepen the understanding of systems governed by discontinuous dynamics.

\section{Model Formulation}
\label{sec:model}
Based on a model introduced in \cite{Campo2017}, we define the impulsive model of competition between two generic species $S_1$ and $S_2$ as:

\begin{subequations}
\begin{align}
&\left.
\begin{cases}\label{eq:equation_1}
        \dfrac{dS_1}{dt}=S_1\left(\psi_1-\dfrac{r_1}{K_1}(S_1+S_2)\right)\left(\dfrac{S_1}{K_0} - 1 \right) - \delta_1 S_1,\\
        \dfrac{dS_2}{dt}=S_2\left(\psi_2-\dfrac{r_2}{K_2}(S_1+S_2)\right)-\delta_2S_2,
\end{cases}
\right.\text{ if } t \neq k\tau,\, k \geq 0\\
&\left.
\begin{cases}\label{eq:equation_2}  
        S_1(t^+)= S_1(t), \\
        S_2(t^+)= S_2(t) + u_k, 
\end{cases}
\right.\text{ if } t = k\tau,\, k \geq 0,
\end{align}
\end{subequations}

\noindent 
with non-negative initial conditions and positive parameters, where $S_1(t)$ and $S_2(t)$ represent the populations of two species competing with each other over time $t$. The parameters $\psi_i$ and $\delta_i$ represent, respectively, the birth and death rates of species $S_1$ and $S_2$ for $i  = 1, 2$, as defined in the original model, while $r_i:= \psi_i - \delta_i$ for $i = 1, 2$ indicates the intrinsic growth rate of both populations. The parameter $K_i$ for $i = 1, 2$ is associated with the carrying capacity of the competing species populations. In this work, we disregard the density dependence of the parameters.

For this impulsive differential equation system, the release period is $\tau$ and $u_k \in U$ denotes the impulsive release of species $S_2$ at time $t = k\tau$. In practice, $u_k$ is limited by the availability of species $S_2$, so the set of possible releases is given by $U:= \{ 0 \leq u_k \leq u_{\text{max}}, \mid k \geq 0 \}$, where $u_{\text{max}} \geq 0$ represents the maximum number of individuals of species $S_2$ that can be released at a given time. The population of $S_1$ immediately after the $k$-th release is given by $t = k\tau^+$, with $S_i(t^+) = \lim_{\epsilon \to 0^+} S_i(t + \epsilon)$, $i=\{1,2\}$.

The system \eqref{eq:equation_1} incorporates the frequency-dependent Allee effect in the first equation, which applies to species $S_1$ \cite{campo2017optimalW10,Campo2017}. This effect is modeled by the critical compensation term $\left(\frac{S_1}{K_0} - 1\right)$, which directly influences the recruitment of individuals of species $S_1$. The term is positive when $S_1(t) > K_0$ and negative when $S_1(t) < K_0$. The parameter $K_0 > 0$, along with $K_1$ (where $0 < K_0 < K_1$), are related to the ``minimum viable population size" (MVPS) commonly observed in models with the Allee effect \cite{barton2011spatialA1, clements2011biologyA2,HASTINGS2013175PopulationDynamics,kanarek2015overcomingA9,stephens1999alleeA7, ufuktepe2022discreteA4}. The MVPS threshold for species $S_1$ is given by $K_b$, while its carrying capacity is indicated by $K_*$. For more details, see \cite{Campo2017}.

For this model, we consider
\begin{equation}\label{eq:equation_3}
    \psi_1 > \delta_1 \mbox{ and } \psi_2 > \delta_2,
\end{equation}
which guarantee a larger number of births than deaths. In addiction, we consider
\begin{equation}\label{eq:equation_4}
    \psi_2 < \psi_1, \, \delta_2 > \delta_1 \text{ and }r_2 < r_1,
\end{equation}
implying that the population of species $S_1$ exhibits greater survival ability than that of species $S_2$.

The system in \eqref{eq:equation_1}, with nonnegative initial conditions, has four steady states, as described in \cite{Campo2017}. In the following, we present these steady states by adapting Theorem $1$ from \cite{Campo2017} to the context of this work.

\begin{theorem}[adapted from Theorem $1$ in \cite{Campo2017}]\label{thm:thm_1}Under the conditions \eqref{eq:equation_4}, the dynamical system \eqref{eq:equation_1} with nonnegative initial conditions has four steady states in the region of biological interest $\mathbb{R}^2_+\setminus\{(0, 0)\}$, namely: 
\begin{itemize}
    \item one nodal repeller $(K_b, 0)$ where
    \begin{equation*}
        K_b = \dfrac{r_1K_0+\psi_1K_1-\sqrt{(r_1K_0+\psi_1K_1)^2-4r_1K_0K_1(\psi_1+\delta_1)}}{2r_1}>0,
    \end{equation*}
    indicates the MVPS threshold for species $S_1$;
    \item one saddle point $(S_1^*, S_2^*)$ of unstable coexistence of both species with coordinates given by
    \begin{align*}
        S_1^* =& \dfrac{K_0\left[\psi_1(K_1-K_2)+\delta_1(K_1+K_2)\right]}{\psi_1(K_1-K_2)+\delta_1K_2}>0,\\
        S_2^* =& K_2-S_1^*>0;
    \end{align*}
    \item two nodal attractors $(0,K_2)$ and $(K_*,0)$, where
    \begin{equation*}
        K_* = \dfrac{r_1K_0+\psi_1K_1+\sqrt{(r_1K_0+\psi_1K_1)^2-4r_1K_0K_1(\psi_1+\delta_1)}}{2r_1}>0,
    \end{equation*}
    defines the carrying capacity of the first species. Only one of these steady states can be reached when $t\to \infty$ according to the initial conditions $S_1(0) > 0$, $S_2(0) > 0$ assigned to the system \eqref{eq:equation_1}, namely:
    \begin{itemize}
        \item[-]If $S_1(0) > K_b$ and $S_2(0) > 0$ then $(K_*,0)$ is
reachable when $t\to \infty$ and the species $S_1$ should persist while the species $S_2$ become extinct.
        \item[-]  If $S_1(0) < K_b$ and $S_2(0) > 0$ then $(0,K_2)$ is
reachable when $t\to \infty$ and the species $S_2$ should persist while the species $S_1$ become extinct.
    \end{itemize}
\end{itemize}
\end{theorem}
In addition, for consistency with the results presented in \cite{{Campo2017}} we consider
\begin{equation*}
    0 < K_0 < K_b < K_2 < K_1 < K_*.
\end{equation*}

Since the system we are working with is a model of interaction between two variables with known terms, we have chosen the model \eqref{eq:equation_1} for this work. When the pulse conditions \eqref{eq:equation_2} are added, it is necessary to study the combined model \eqref{eq:equation_1}--\eqref{eq:equation_2} to ensure that, even with the discontinuity, the system remains well-posed.

In the next section, we will explore the dynamic behavior of the impulsive differential equations system, including an analysis of the local and global stability of the free solution of species $S_1$. We will also present numerical simulations that illustrate these behaviors and the corresponding results, contributing to a deeper understanding of the dynamics between the species.

\section{Model Dynamics Analysis}
\label{sec:dynamics}
In this section, we begin by presenting some essential mathematical de\-finitions and tools for analyzing the impulsive differential equations system. These tools are fundamental for establishing results related to positivity, existence, uniqueness, and boundedness of solutions. We perform a detailed analysis of the system, focusing on the results obtained and their implications for the behavior of solutions over time.

Additionally, when considering the release of individuals from species $S_2$ to replace the population of $S_1$. We prove the existence and uniqueness of the solution $(0, \overline{S}_2)$ and analyze the conditions for its stability. The behavior of this solution is of great importance, as it allows us to understand how species $S_2$ behaves in the absence of direct interactions with other species. Through this investigation, we aim to better understand the conditions that ensure the maintenance of the free solution, which may have significant implications for the system's dynamics and species interactions.

\subsection{Fundamental Concepts}

In this subsection, we present some definitions and tools that will be used to analyze the existence and stability of impulsive periodic solutions of system \eqref{eq:equation_1}–\eqref{eq:equation_2}. Throughout the paper, we adopt the notations $\mathbb{R}_+=[0, \infty)$, $\mathbb{R}_+^2 = \{x = (x_1, x_2) \in \mathbb{R}^2 : x_1 \geq 0, \ x_2 \geq 0\}$, and denote by $g = (g_1, g_2)^T$ the vector field on the right-hand side of \eqref{eq:equation_1}.

\begin{definition}{\rm (in \cite{Laksh1989})}\label{def:1}
      Let $V : \mathbb{R}_+ \times \mathbb{R}_+^2 \to \mathbb{R}_+$. We say that $V$ belongs to the class $\mathcal{V}_0$ if it is continuous on $(k\tau, (k + 1)\tau] \times \mathbb{R}_+^2$ for all $k \in \mathbb{N}$, and if the following limit exists and is finite for every $x \in \mathbb{R}_+^2$:
\begin{equation*}
    \lim_{(t, y) \to (k\tau^+, x)} V(t, y) = V(k\tau^+, x).
\end{equation*}

\end{definition}

\begin{definition}{\rm (in \cite{Laksh1989})}\label{def:2} 
Let $V \in \mathcal{V}_0$. For $(t, x) \in (k\tau, (k+1)\tau] \times \mathbb{R}_+^2$, the upper right-hand derivative of $V(t,x)$ with respect to the impulsive system \eqref{eq:equation_1}–\eqref{eq:equation_2} is defined by
\begin{equation*}
    D^+ V(t,x) = \lim_{h \to 0} \sup \frac{1}{h}\left[V(t+h,x+h g(t,x)) - V(t,x)\right].
\end{equation*}
\end{definition}

\begin{definition}{\rm (in \cite{Laksh1989})}\label{def:3}
Let $\varrho(t) = \varrho(t, t_0, x_0)$ be a solution of system \eqref{eq:equation_1}–\eqref{eq:equation_2} defined on an interval $[t_0, t_0 + \ell)$. The function $\varrho(t)$ is called a \emph{maximal solution} if, for any other solution $x(t, t_0, x_0)$ defined on the same interval, we have
    $$x(t) \leq \varrho(t), \quad \text{for all } t \in [t_0, t_0 + \ell).$$
\end{definition}
 A \emph{minimal solution} $\rho(t)$ is defined analogously by reversing the inequality in Definition \ref{def:3}.

Using Definitions \ref{def:1}, \ref{def:2}, and \ref{def:3}, we can now state a comparison theorem for impulsive differential equations.

\begin{theorem}\label{thm:thm_2} {\rm (}Comparison theorem {\rm in \cite{Laksh1989})}:
Let $m \in \mathcal{V}_0$, and suppose that
\begin{align*}
    D^+ m(t) &\leq v(t, m(t)), && t \neq t_k, \quad k = 1,2,\dots \\
    m(t_k^+) &\leq \varphi_k(m(t_k)), && t = t_k, \quad k = 1,2,\dots \nonumber
\end{align*}
where each $\varphi_k \in \mathcal{C}(\mathbb{R}, \mathbb{R})$ is non-decreasing.

Let $\varrho(t)$ be the maximal solution of the scalar impulsive system
\begin{align}\label{eq:equation_5}
    \dot{u}(t) &= v(t,u), && t \neq t_k \\
    u(t_k^+) &= \varphi_k(u(t_k)), && t = t_k, \quad t_k > t_0 \geq 0 \nonumber \\
    u(t_0) &= u_0, \nonumber
\end{align}
which exists on $[t_0, \infty)$. Then, if $m(t_0^+) \leq u_0$, it follows that $m(t) \leq \varrho(t)$ for all $t \geq t_0$.

A symmetric result holds if all inequalities are reversed and the maps $\varphi_k$ are non-increasing.

\end{theorem}

\begin{remark}
    Although Definition~\ref{def:1} introduces the class $\mathcal{V}_0$ for functions $V(t, x)$ that depend on both time and state variables, scalar functions of time, such as $m(t)$ in Theorem~\ref{thm:thm_2}, are also said to belong to $\mathcal{V}_0$ since they satisfy the corresponding regularity properties at impulsive moments. More precisely, a scalar function $m : \mathbb{R}_+ \to \mathbb{R}$ belongs to class $\mathcal{V}_0$ if it satisfies:
    \begin{itemize}
        \item $m(t)$ is continuous on each interval $(t_k, t_{k+1}]$;
        \item The right-hand limit $\displaystyle \lim_{t \to t_k^+} m(t)$ exists and is finite for all $k$.
    \end{itemize}
    These conditions ensure that $m$ behaves regularly at the discontinuities induced by the impulses, analogously to the original definition of class $\mathcal{V}_0$.
\end{remark}

Note that if the function $v$ in Theorem \ref{thm:thm_2} is smooth enough to ensure existence and uniqueness of solutions for the initial value problem \eqref{eq:equation_5}, then $\varrho(t)$ is the unique solution.

\subsection{Behavior of system solutions}\label{subs:behavior}
As we are working with population dynamics, we must make sure that the solutions of system \eqref{eq:equation_1}-\eqref{eq:equation_2} are non-negative.
\begin{proposition}\label{prop:prop_1}
    Let $(S_1(0),S_2(0))$ be a non-negative initial condition, and let $(S_1(t),S_2(t))$ be a solution to the system \eqref{eq:equation_1}-\eqref{eq:equation_2}. Then, $(S_1(t),S_2(t))$ remains non-negative for all $t\geq 0$.
\end{proposition}
\begin{proof}
     Note that $\dfrac{dS_1}{dt} = 0$ whenever $S_1(t) = 0$. This implies that the solution cannot cross the $S_1$-axis. Therefore, if $S_1(0) \geq 0$, it follows that $S_1(t) \geq 0$ for all $t \geq 0$. Similarly, if $S_2(0) \geq 0$, then $S_2(t) \geq 0$ for all $t \geq 0$.
\end{proof}

In the following proposition, the smoothness of the right-hand side of the system \eqref{eq:equation_1}--\eqref{eq:equation_2}, combined with Definition~\ref{def:1}, ensures the existence and uniqueness of solutions to the model.

\begin{proposition}\label{prop:prop_2} For each non-negative initial condition and each release a\-mount  $u_k \in U$ of individuals from species $S_2$, the system \eqref{eq:equation_1}–\eqref{eq:equation_2} has a unique solution defined on the interval $[0, \infty)$.
\end{proposition}

\begin{proof}
Suppose $(S_1(t), S_2(t))$ is a solution to the system \eqref{eq:equation_1}--\eqref{eq:equation_2}. It is continuous on the intervals $(k\tau, (k+1)\tau]$ for $k \geq 0$, indicating that it remains continuous between each pair of pulses. Furthermore, there exist limits defined as follows:
\begin{equation*}
S_1(k\tau^+) = \lim_{\epsilon \to 0^+} S_1(k\tau + \epsilon), \quad S_2(k\tau^+) = \lim_{\epsilon \to 0^+} S_2(k\tau + \epsilon).
\end{equation*}
Consequently, the existence and uniqueness of these solutions are guaranteed by the smoothness of the functions:
\begin{align*}
&g_1(S_1, S_2) = S_1\left(\psi_1 - \frac{r_1}{K_1}(S_1 + S_2)\right)\left(\frac{S_1}{K_0} - 1\right) - \delta_1 S_1,\\
&g_2(S_1, S_2) = S_2\left(\psi_2 - \frac{r_2}{K_2}(S_1 + S_2)\right) - \delta_2 S_2.
\end{align*}
\end{proof}

Due to the biological context of system \eqref{eq:equation_1}--\eqref{eq:equation_2}, it is essential to ensure that the solutions remain bounded over time, as unbounded growth or negative population sizes would not be biologically meaningful. The model captures key ecological mechanisms such as intra- and interspecific competition, mortality, and impulsive interventions (e.g., periodic releases). To better understand the dynamics and provide a foundation for proving boundedness, we first consider an auxiliary system that describes the evolution of species $S_2$ in the absence of individuals of species $S_1$.
\begin{align}
\begin{cases}\label{eq:equation_6}
\dfrac{dZ_2}{dt}(t) = Z_2\left(\psi_2 - \dfrac{r_2}{K_2} Z_2\right) - \delta_2 Z_2.
\quad t \neq k\tau, \ k \geq 0  \\    
Z_2(t^+)=Z_2(t) + u_k,  \qquad t = k\tau, \ u_k \in U\\
Z_2(0^+)=Z_2(0).
\end{cases}
\end{align}

This simplified system captures the autonomous dynamics of $S_2$ under impulsive control, without interaction terms involving $S_1$. We will show that, for each initial condition $Z_2(0)$ such that $Z_2(0) = S_2(0)$ and for each set of releases $U$, there exists a unique $\tau$-periodic solution that is globally asymptotically stable. This result is a key step toward establishing uniform bounds on the full system's solutions, and it ensures the existence and global stability of the $S_1$-free solution of system \eqref{eq:equation_1}--\eqref{eq:equation_2}.
\begin{theorem}\label{thm:thm_3}Given the auxiliary system \eqref{eq:equation_6}, for each $Z_2(0) \geq 0$ and each $u_k \in U$, there exists a unique positive $\tau$-periodic solution $\overline{Z}_2(t)$, expressed by:
\begin{equation}\label{eq:equation_7}
    \overline{Z}_2(t) = \dfrac{K_2 Z_2^{+} e^{r_2 (t-k\tau)}}{Z_2^{+}\left(e^{r_2 (t-k\tau)} -1\right) + K_2}, \quad k\tau < t \leq (k+1)\tau, \, k \geq 0,
\end{equation}
where
\begin{equation}\label{eq:equation_8}
   Z_2^{+} = \dfrac{1}{2}\left[(u_k+K_2)+\sqrt{(u_k+K_2)^2+4\dfrac{u_kK_2}{e^{r_2 \tau}-1}}\right], \,k \geq 0.
\end{equation}
Furthermore, the solution $\overline{Z}_2(t)$ is globally asymptotically stable.
\end{theorem}
The proof of Theorem \ref{thm:thm_3} is given in the \eqref{appendix:A}.

\begin{corollary}\label{cor:cor_1}
    Let $\overline{Z}_2(t)$ be the $\tau$-periodic solution of system \eqref{eq:equation_6}. Then for every solution $Z_2(t)$ of problem \eqref{eq:equation_6}, 
    \begin{equation*}
        Z_2(t) \to \overline{Z}_2(t), \mbox{ as } t \to \infty,
    \end{equation*}where $\overline{Z}_2(t)$ is given by \eqref{eq:equation_7}.
\end{corollary}
\begin{proof}
    For every solution $Z_2(t)$ of the system \eqref{eq:equation_6}, we have $Z_2(t) \to \overline{Z}_2(t)$ as  $t \to \infty$  follows directly from the \Cref{thm:thm_3} that establishes the global asymptotic stability of  $\overline{Z}_2(t)$.
\end{proof}

Now, we are interested in proving the uniform boundedness of the solutions of system \eqref{eq:equation_1}-\eqref{eq:equation_2}. To this end, we first establish a lemma that provides an upper bound for $S_1(t)$, based on its relationship with the carrying capacity $K_*$, as discussed in \cite{Campo2017}. Then, we employ the auxiliary system \eqref{eq:equation_6}, along with Theorem~\ref{thm:thm_3} and its Corollary~\ref{cor:cor_1}, to derive an upper bound for $S_2(t)$.

From a biological perspective, proving that the population sizes remain bounded ensures that the model describes realistic population dynamics under impulsive interventions. From a mathematical point of view, uniform boundedness plays a key role in establishing the global asymptotic stabi\-lity of the $S_1$-free solution and is also a crucial hypothesis for proving the existence of an optimal control, which will be addressed in the next section.

\begin{lemma}\label{lem:lem_1}
     Let $(S_1(t), S_2(t))$ be a solution of system \eqref{eq:equation_1}--\eqref{eq:equation_2}, with positive parameters, $u_k \in U$, and non-negative initial conditions. Then, there exists a positive constant $M_1$ such that $S_1(t) \leq M_1$ for all $t \geq 0$.
\end{lemma}
The proof of Lemma \ref{lem:lem_1} is given in the \eqref{appendix:B}.

\begin{theorem}\label{prop:prop_3}
    Let $(S_1(t),S_2(t))$ be a solution of system \eqref{eq:equation_1}-\eqref{eq:equation_2}, with positive parameters, $u_k \in U$ and non-negative initial conditions. Then $(S_1(t),S_2(t))$ is uniformly  bounded.
\end{theorem}

The proof of Theorem \ref{prop:prop_3} is given in the \eqref{appendix:C}.\\

With the support of Theorem \ref{thm:thm_3}, we establish the existence of a unique $\tau$-periodic solution of system \eqref{eq:equation_1}–\eqref{eq:equation_2} corresponding to the absence of species $S_1$.

\begin{theorem}\label{thm:thm_4} 
Let $(S_1(0), S_2(0))$ be non-negative initial conditions. For some choice of $u_k \in U$, the pair $(0, \overline{S}_2(t))$ is the unique positive, $\tau$-periodic, $S_1$-free solution of system \eqref{eq:equation_1}--\eqref{eq:equation_2}.
\end{theorem}

\begin{proof}
 Note that when $S_1(t) = 0$, the second equation of system \eqref{eq:equation_1} reduces to
$\dfrac{dS_2}{dt} = \dfrac{dZ_2}{dt}$, where $\dfrac{dZ_2}{dt}$ corresponds to the first equation of the auxiliary system \eqref{eq:equation_6}. The solution to this equation is given explicitly by \eqref{eq:equation_7}.

Therefore, the dynamics of $S_2(t)$ under the condition $S_1(t) = 0$ coincide exactly with those of $Z_2(t)$. By applying Theorem \ref{thm:thm_3}, which guarantees the existence and uniqueness of a positive $\tau$-periodic solution for system \eqref{eq:equation_6}, we conclude that there exists $u_k \in U$ such that $(0, \overline{S}_2(t))$ is the unique $\tau$-periodic, $S_1$-free solution of system \eqref{eq:equation_1}--\eqref{eq:equation_2}, where $\overline{S}_2(t) = \overline{Z}_2(t)$ for all $t \geq 0$.
\end{proof}

\subsection{Stability of $S_1$-free periodic solution}
\label{sec:stability}
 Due to the Allee effect present in the equation governing species $S_1$, we observed in the previous section that when $S_1(0) < K_b$, the population of $S_1$ tends to go extinct in the absence of external intervention. Therefore, it is crucial to analyze the stability of the solution $(0, \overline{S_2}(t))$ when $S_1(0) > K_b$. Understanding the behavior of this solution is essential for evaluating the viability of species $S_2$ in the environment. In this subsection, we investigate the global asymptotic stability of the periodic solution with no individuals of species $S_1$ in the system \eqref{eq:equation_1}-\eqref{eq:equation_2}. To analyze this, we use a comparison argument with auxiliary systems whose dynamics are known. Additionally, we rely on the fact that all solutions of the system are uniformly bounded to ensure that the convergence is global.

\begin{theorem}\label{thm:thm_6}
The $S_1$-free periodic solution $(0, \overline{S}_2(t))$ of the system \eqref{eq:equation_1}–\eqref{eq:equation_2} is globally asymptotically stable if
\begin{equation}\label{eq:equation_9}
    \overline{S}_2(t) > K_1\quad \text{for all } t \geq 0.
\end{equation}
\end{theorem}
\begin{proof}
From the second equation of the model, we obtain the inequality:
\begin{equation*}
 \dfrac{dS_2}{dt}\leq\left(\psi_2-\dfrac{r_2}{K_2}S_2\right)S_2 - \delta_2 S_2.   
\end{equation*}
Therefore, we can apply the auxiliary system~\eqref{eq:equation_6} for comparison. From Corollary~\ref{cor:cor_1}, we know that
\begin{equation*}
    \lim_{t \to \infty} Z_2(t) = \overline{Z}_2(t).
\end{equation*}
Thus, for any $\epsilon > 0$ sufficiently small, there exists $t_1 > 0$ such that 
\begin{equation*}
    Z_2(t)<\overline{Z}_2(t)+\epsilon \quad \text{for all } t > t_1.
\end{equation*}
It follows from the Comparison Theorem~\ref{thm:thm_2} that $S_2(t) \leq Z_2(t)$. Then, since $Z_2(t) < \overline{Z}_2(t) + \epsilon$ and $S_2(0) = Z_2(0)$, we conclude that
\begin{equation}\label{eq:equation_10}
    S_2(t) \leq \overline{Z}_2(t) + \epsilon, \quad \text{for all } t > t_1.
\end{equation}

Now, suppose that $S_1(0) > K_b > K_0$. In this case, $\left(\frac{S_1}{K_0} - 1\right) > 0$. Consider the inequality: 
\begin{align*}
    \dfrac{dS_1}{dt} &\leq S_1\left[\left(\psi_1 - \dfrac{r_1}{K_1}S_2\right)\left(\dfrac{S_1}{K_0} - 1\right) - \delta_1\right].
\end{align*}To ensure the stability of the equilibrium solution $S_1(t) = 0$, we require that
% Thus, for $S_1(t) \to 0$ as $t \to \infty$, it is sufficient that
\begin{equation}\label{eq:equation_11}
    \left(\psi_1 - \dfrac{r_1}{K_1}S_2\right)\left(\dfrac{S_1}{K_0} - 1\right)< \delta_1.
\end{equation}
Given \eqref{eq:equation_10}, the inequality~\eqref{eq:equation_11} is guaranteed if
\begin{equation*}
    \psi_1 - \dfrac{r_1}{K_1}(\overline{Z}_2(t)+\epsilon) < \dfrac{\delta_1}{\left(\dfrac{S_1}{K_0} - 1\right)}.
\end{equation*}
Since $S_1(0) >K_b> K_0$, we have $\left(\dfrac{S_1}{K_0} - 1\right) > 1$ for sufficiently large values of $S_1$, which implies that
\begin{equation*}
 \dfrac{\delta_1}{\left(\dfrac{S_1}{K_0} - 1\right)} < \delta_1.   
\end{equation*}
Hence, the inequality
\begin{align*}
     \psi_1 - \dfrac{r_1}{K_1}(\overline{Z}_2(t)+\epsilon)  < \delta_1     
\end{align*}
is also satisfied. That is, we can rewrite the condition as
\begin{equation*}
      \overline{Z}_2(t) > K_1\quad \text{for all } t \geq 0,
\end{equation*}
since by the model hypothesis $r_1 = \psi_1 - \delta_1$.

Substituting~\eqref{eq:equation_9} into the first equation of system~\eqref{eq:equation_1}, we obtain
\begin{align*}
    \dfrac{dS_1}{dt} \leq S_1\left[\left(\psi_1 - \dfrac{r_1}{K_1}(S_1+K_1)\right)\left(\dfrac{S_1}{K_0} - 1\right) - \delta_1\right].
\end{align*}
We now consider the following comparison system:
\begin{align*}
\begin{cases}
    \dfrac{dZ_1}{dt} = Z_1\left[\left(\psi_1 - \dfrac{r_1}{K_1}(Z_1+K_1)\right)\left(\dfrac{Z_1}{K_0} - 1\right) - \delta_1\right],\\
    Z_1(0) = S_1(0),
\end{cases}
\end{align*}
for which we know that $Z_1(t) \to 0$ as $t \to \infty$, by construction. Therefore, by the Comparison Theorem, if condition~\eqref{eq:equation_9} holds, then for $\epsilon > 0$ sufficiently small, there exists $t_2 > t_1$ such that 
\begin{equation}\label{eq:equation_12}
    S_1(t) \leq Z_1(t) \leq \epsilon, \quad \text{for all } t > t_2.
\end{equation}

Substituting \eqref{eq:equation_12} into the second equation of system~\eqref{eq:equation_1}, we obtain:
\begin{align*}
\begin{cases}
\dfrac{dS_2}{dt}(t)\geq S_2\left(\psi_2-\dfrac{r_2}{K_2}(\epsilon+S_2)\right)-\delta_2 S_2,
\quad t \neq k\tau, \ k \geq 0, \\    
S_2(t^+)=S_2(t) + u_k,  \qquad t = k\tau, \ u_k \in U,\\
S_2(0^+)=S_2(0).
\end{cases}
\end{align*}
Due to the continuity of the right-hand side, and by an argument analogous to that used in Theorem~\ref{thm:thm_3}, we conclude that for $\epsilon > 0$ sufficiently small, there exists $t_3 > t_2$ such that 
\begin{equation*}
    \overline{Z}_2(t)-\epsilon \leq S_2(t), \quad \text{for all } t > t_3.
\end{equation*}

Finally, if condition~\eqref{eq:equation_9} is satisfied, then for $\epsilon > 0$ we have:
\begin{equation*}\
    0 \leq S_1(t) \leq \epsilon \quad \text{and} \quad \overline{Z}_2(t)-\epsilon \leq S_2(t) \leq \overline{Z}_2(t)+\epsilon, \quad \text{for all } t > t_3.
\end{equation*}
Letting $\epsilon \to 0$, we conclude that
\begin{equation*}
    S_1(t) \to 0 \quad \text{and} \quad S_2(t) \to \overline{Z}_2(t):= \overline{S}_2(t)  \quad \text{as } t \to \infty.
\end{equation*}

This implies that $(S_1(t), S_2(t)) \to (0, \overline{S}_2(t))$ as $t \to \infty$ for any admissible initial condition satisfying $S_1(0) > K_0$, provided that condition~\eqref{eq:equation_9} holds. Furthermore, since the solutions remain uniformly bounded and converge regardless of the initial condition (within the biologically relevant domain), the solution $(0, \overline{S}_2(t))$ is globally attractive and stable in the sense of Lyapunov. Therefore, it is globally asymptotically stable.
\end{proof}

\subsection{A method to select $u_k$ (a sufficient condition)}\label{subs:subsec3.4}
In this subsection, we make use of the global stability condition \eqref{eq:equation_9} to determine a sufficient value for the number of individuals of species $S_2$ to be released, denoted by $u_k$, to stabilize the $S_2$ population and drive the $S_1$ population to extinction.

Let
\begin{align*}
    \overline{S}_2(t) = \dfrac{K_2 Z_{2}^{+} e^{r_2 (t-k\tau)}}{Z_{2}^{+}\left(e^{r_2 (t-k\tau)} - 1\right) + K_2}, \quad  k\tau < t \leq (k+1)\tau, \quad k \geq 0,
\end{align*}
where $Z_2^+$ depends on $u_k$ and is given by \eqref{eq:equation_8}. The maximum value of $\overline{S}_2(t)$ over each interval $(k\tau, (k+1)\tau]$ occurs at $t = (k+1)\tau$, yielding
\begin{equation*}\label{eq:equation_41}
    \overline{S}_2^{\max} = \dfrac{K_2 Z_{2}^{+} e^{r_2 \tau}}{Z_{2}^{+}\left(e^{r_2\tau} -1\right)+K_2},
\end{equation*}
so that
\begin{equation*}
    \overline{S}_2(t) \leq \overline{S}_2^{\max}, \quad \text{for all } t \geq 0.
\end{equation*}

In order to solution $(0, \overline{S}_2(t))$ be globally asymptotically stable, it is sufficient that $\overline{S}_2(t) > K_1$. Thus, requiring
\begin{equation*}\label{eq:equation_43}
    \overline{S}_2^{\max} > K_1,
\end{equation*}
and performing algebraic manipulations, we obtain the following sufficient condition to guarantee the stabilization of $S_2$ and the extinction of $S_1$:
\begin{equation*}
    u_k > \eta(\tau), \quad \forall \tau \geq 0,
\end{equation*}
where
\begin{equation}\label{eq:equation_13}
    \eta(\tau) = \frac{(e^{r_2 \tau}-1)\phi(\tau)(\phi(\tau)-K_2)}{K_2 + \phi(\tau)(e^{r_2 \tau}-1)} \quad \text{and} \quad \phi(\tau) = \frac{{K}_1 K_2}{e^{r_2 \tau} K_2 - {K}_1(e^{r_2 \tau}-1)}.
\end{equation}

Observing $\eta(\tau)$ may assume negative values for certain values of $\tau > 0$ and since the biological context requires $u_k > 0$, we restrict our attention to values of $\tau$ for which $\eta(\tau)$ is positive, and define a conservative lower bound:
\begin{equation}\label{eq:equation_14}
    u_k > \max_{\tau > 0} \eta(\tau) > 0,
\end{equation}
to ensure that the sufficient condition holds uniformly. We now demonstrate that this maximum exists.

\begin{proposition}\label{prop:prop_4}
Let $\eta(\tau)$ and $\phi(\tau)$ be defined for $\tau \geq 0$ as in \eqref{eq:equation_13}. Then, the function $\eta(\tau)$ attains a global maximum on $[0, \infty)$.
\end{proposition}

\begin{proof}
The function $\eta(\tau)$ is continuous and differentiable on the interval $[0, \infty)$. Moreover, we have:
$$
\lim_{\tau \to 0} \phi(\tau) = K_1, \quad \lim_{\tau \to 0} \eta(\tau) = 0,
\quad \text{and} \quad \lim_{\tau \to \infty} \phi(\tau) = 0, \quad \lim_{\tau \to \infty} \eta(\tau) = 0.
$$
Although the interval is unbounded above, an analysis of the function $\eta(\tau)$ shows that it reaches a local (and global) maximum at some $\tau = \tau_{\text{max}} \in (0, \infty)$, which completes the proof.
\end{proof}

The analysis carried out in this section revealed the conditions for the global stability of the $S_1$-free solution, providing a solid theoretical foundation for the dynamic behavior of the impulsive model. In the next section, we address the optimal control problem, where the previously introduced model is used as the control system, and the results obtained throughout this section are fundamental for ensuring the existence of an optimal solution. This problem is formulated to determine strategies that maximize the effectiveness of population dynamics control. The numerical results of these approaches, as well as their implications, will be presented in \Cref{sec:numericalresults}.

\section{Impulsive Optimal Control Problem}
In this section, we formulate an optimal control problem involving the impulsive release of individuals of species $S_2$. The control actions consist of fixed-frequency releases, and the objective is to determine the optimal number of individuals to be introduced at each intervention time. A central question guides this formulation: What is the minimum number of $S_2$ individuals that must be released, at a fixed frequency, to guarantee their fixation in the target population while minimizing the overall intervention cost?

This control strategy is applied over a finite time horizon, denoted by $[0, T]$, where $T$ represents the final time of intervention. Within this interval, up to $N$ releases may occur, each corresponding to an element of the admissible control set, defined as
$\bar{U}:=\{u_k \in \mathbb{R} \mid 0 \leq u_k \leq u_{max}, \, k = 1,2,..., N\}$, where $u_{\max}$ satisfies $u_{\max} \geq \max\limits_{\tau>0} \eta(\tau)$ with $\eta(\tau)$ specified in the Subsection \ref{subs:subsec3.4}. The release period $\tau$ is considered fixed, determined by operational constraints such as logistical planning, environmental conditions, release policies, and economic factors, which limit the flexibility to vary the timing between consecutive releases.

Given this framework, the optimal control problem can be stated as follows: Find an optimal control $u^* = (u^*_{k})_{k=1}^{N}$, with each $u_k \in \bar{U}$, that minimizes the cost functional

\begin{small}
\begin{equation}\label{eq:equation_15}
    J(u) = C \sum_{k=1}^{N} u_k,
\end{equation}
subject to
\begin{align}
&\left\{
\begin{aligned}
    \dfrac{dS_1}{dt} &= S_1\left(\psi_1 - \dfrac{r_1}{K_1}(S_1 + S_2)\right)\left(\dfrac{S_1}{K_0} - 1 \right) - \delta_1 S_1, \\
    \dfrac{dS_2}{dt} &= S_2\left(\psi_2 - \dfrac{r_2}{K_2}(S_1 + S_2)\right) - \delta_2 S_2,
\end{aligned}
\right. && \text{if } t \neq k\tau,\; k = 1, 2, \dots, N, \label{eq:equation_16} \\
&\left\{
\begin{aligned}
    S_1(t^+) &= S_1(t), \\
    S_2(t^+) &= S_2(t) + u_k,
\end{aligned}
\right. && \text{if } t = k\tau,\; k = 1, 2, \dots, N, \label{eq:equation_17} \\
&\left.\begin{cases}\label{eq:equation_18}
    S_1(T) < K_b
\end{cases}\right.
\end{align}
with initial conditions:
\begin{equation*}
    S_1(0) \geq 0 \quad \text{and} \quad S_2(0) \geq 0.
\end{equation*}
\end{small}

The condition $u_{max} \geq \max\limits_{\tau>0} \eta(\tau)$ is established to ensure that the admissible control set $\bar{U}$ is non-empty, which in turn guarantees the existence of at least one feasible value of $u_k$ satisfying the problem described in equations \eqref{eq:equation_15}--\eqref{eq:equation_18}, as will be shown in Proposition \ref{prop:prop_5}.

We consider the following aspects of the control problem under investigation:
\begin{itemize}
    \item[i.] The objective function $J(u)$ accounts for the sum of all releases, with the constant $C$ representing the cost of intervention;
    \item[ii.] $J(u)$ is subject to the model dynamics defined by equations \eqref{eq:equation_16}--\eqref{eq:equation_17}, which characterize the control system;
    \item[iii.] The constraint \eqref{eq:equation_18} aims to reduce the $S_1$ population at the final time \( T \) below the threshold \( K_b \), associated with the Allee effect (see \cite{Campo2017,barton2011spatialA1}), to ensure the fixation of the $S_2$ species.
\end{itemize}

Next, we present a proposition that demonstrates that $ \bar{U}$ is a non-empty and compact set that satisfies the final constraint of the problem. This property is essential to ensure the existence of a viable solution within the established conditions.

\begin{proposition}\label{prop:prop_5}
    The set of admissible controls $ \bar{U} $ is non-empty and compact.
\end{proposition}
\begin{proof}
In Section~\ref{subs:behavior}, we showed that the system \eqref{eq:equation_16}--\eqref{eq:equation_17} is well-posed for any $ u_k \in U $. Since $ \bar{U} \subset U$, the system is also well-posed for any $ u_k \in \bar{U} $, ensuring the existence of a unique solution.

By hypothesis, $u_{max} \geq \max\limits_{\tau>0} \eta(\tau)$ is given in equation \eqref{eq:equation_13}. Choosing $ k^* \in \{1, \dots, N\} $ such that $ u_{k^*} =  \max\limits_{\tau>0} \eta(\tau)$, we obtain $ \overline{S}_2 > K_1 $ for this release, as shown in Theorem~\ref{thm:thm_6}. Under this condition, the dynamics of $ S_1 $,
$$
\dfrac{dS_1}{dt} = S_1\left(\psi_1 - \dfrac{r_1}{K_1}(S_1 + S_2)\right)\left(\dfrac{S_1}{K_0} - 1 \right) - \delta_1 S_1,
$$
lead to a decline of $ S_1 $ toward zero in finite time, satisfying the constraint $ S_1(T) < K_b $.

Moreover, since $ u_k \in [0, u_{\max}] $, the set $ \bar{U} $ is closed and bounded in $ \mathbb{R}$, and therefore compact. We conclude that $ \bar{U} $ is non-empty and compact, containing at least one control $ u_{k^*} \in \bar{U} $ that satisfies the system and the terminal constraint~\eqref{eq:equation_18}.
\end{proof}

Based on Theorem 5.1 in Section III of \cite{berkovitz2013optimal}, we state the following result, which guarantees the existence of an optimal control for the problem under consideration.
\begin{theorem}{\rm (Existence of optimal control)}\label{thm:thm_7}
Consider the control problem defined by equations \eqref{eq:equation_15}-\eqref{eq:equation_18}, with the admissible control set $\bar{U} \subset \mathbb{R}$.

If the following conditions are satisfied:
\begin{enumerate}
    \item The admissible control set $\bar{U}$ is nonempty and compact;
    \item For each control $u = (u_k)_{k=1}^N \in \bar{U}$, the impulsive system \eqref{eq:equation_16}-\eqref{eq:equation_17} has a unique, positive, and uniformly bounded solution;
    \item The cost functional $J(u) = C \sum_{k=1}^{N} u_k$ is continuous with respect to $u$,
\end{enumerate}
then there exists an optimal control $u^* = (u_k^*)_{k=1}^N \in \bar{U}$ that minimizes the cost functional $J(u)$ subject to the system dynamics and constraints.
\end{theorem}

\begin{proof}
    The first condition has already been verified in Proposition~\ref{prop:prop_5}. In Subsection~\ref{subs:behavior}, we established the existence, uniqueness, positivity, and uniform boundedness of the solutions to the system \eqref{eq:equation_16}--\eqref{eq:equation_17} for each $u = (u_k)_{k=1}^N \in \bar{U} \subset U$, which satisfies the second condition. 

    Furthermore, the cost functional $J(u) = C \sum_{k=1}^{N} u_k$ is a linear function of the control variables $u_k$, and since $\bar{U}$ is compact, it follows that $J$ is continuous on $\bar{U}$, satisfying condition (iii). 

    Since $J$ is continuous and $\bar{U}$ is nonempty and compact, we can apply the Weierstrass Theorem, which ensures that every continuous function defined on a compact set attains its minimum on that set. Therefore, there exists $u^* = (u_k^*)_{k=1}^N \in \bar{U}$ such that
    \[
    J(u^*) = \min_{u \in \bar{U}} J(u),
    \]
    concluding that $u^*$ is an optimal control for the problem defined by equations \eqref{eq:equation_15}-\eqref{eq:equation_18}.
\end{proof}

Based on these results, for any initial configuration of the system, it is possible to find an optimal control policy that satisfies the necessary conditions, thereby ensuring the viability of the proposed solution. The existence of an optimal control is a central topic in control theory, relying on classical assumptions such as the compactness of the admissible control set and the continuity of the system dynamics (see \cite{berkovitz2013optimal, cesari1983optimization, fleming1975deterministic}). These assumptions are fulfilled in our setting, as demonstrated by the formulation of the problem and the compactness of the admissible set $\bar{U}$.

In the next section, we present numerical results that illustrate the efficiency of the proposed control strategy.

\section{Numerical Results}\label{sec:numericalresults}
This section presents the numerical results obtained to validate the theoretical findings and to investigate the behavior of the optimal control problem through numerical simulations. First, we assess the consistency of the theoretical results by comparing them with numerical approximations. Next, we explore the numerical solution of the optimal control problem, illustrating how the proposed approach performs under different parameter configurations and highlighting key observations. All simulations use the parameters detailed in \Cref{tab:tab_1} from \cite{Campo2017}.

\begin{table}[htpb]
\centering
{\footnotesize
\caption{Parameters for the model \eqref{eq:equation_1}-\eqref{eq:equation_2}}\label{tab:tab_1}
\begin{tabular}{lccc}
\toprule
\textbf{Parameter} & \textbf{Value} & \textbf{Range} & \textbf{Description} \\
\midrule
$\psi_1$ & 0.32667 & 0.28 - 0.38 & Birth rate of species $S_1$ \\
$\psi_2$ & 0.21333 & 0.18 - 0.25 & Birth rate of species $S_2$ \\
$\delta_1$ & 0.03333 & $1/8$ - $1/42$ & Death rate of species $S_1$ \\
$\delta_2$ & 0.06666 & $2/8$ - $2/42$ & Death rate of species $S_2$ \\
$K_1$ & 374 & - & Related to the carrying capacity of $S_1$ \\
$K_2$ & 300 & - & Carrying capacity of species $S_2$ \\
$K_0$ & 30 & - & Threshold population for species interaction \\
\bottomrule
\end{tabular}}
\end{table}

For the simulations, we illustrate the species in question, with $S_1$ representing wild female Aedes aegypti mosquitoes and $S_2$ representing female mosquitoes infected with Wolbachia bacteria \cite{almeida2019mosquitoW5,cardona2020wolbachiaW7,zara2016estrategiasW2}. The model without impulsive control, proposed by \cite{Campo2017}, was applied in the context of Wolbachia-based Aedes aegypti control. Here, we extend this model to address any two species or subspecies, provided that their populations satisfy conditions \eqref{eq:equation_3} and \eqref{eq:equation_4}. As an example, we return to the specific case involving Wolbachia.

\subsection{Numerical consistency of theoretical results}

To ensure a comprehensive assessment of the dynamics between species $S_1$ and $S_2$ (wild and Wolbachia-infected female Aedes aegypti mosquitoes), we performed simulations using four different initial conditions. The specific objective of this subsection is to demonstrate the global stability of the solution $(0, \overline{S}_2(t))$, based on the condition established in \Cref{thm:thm_6}.

The simulations were conducted over 180 days, using the parameters listed in \Cref{tab:tab_1}, with Python and its libraries ensuring precision and efficiency in solving the system of equations. We implemented the fourth-order Runge–Kutta method with the NumPy library, including adaptations to incorporate the impulsive jumps, providing a robust numerical solution.

\begin{figure}[!ht]
    \centering
    % Primeira imagem com etiqueta (a)
    \begin{minipage}{0.48\textwidth}
        \centering
        \includegraphics[width=\textwidth]{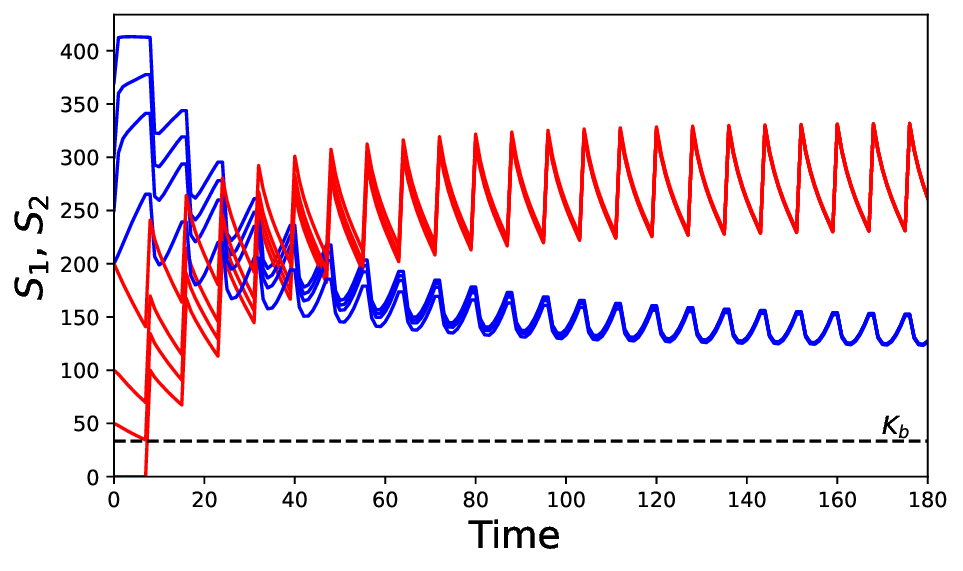}
        \caption*{(a)}
    \end{minipage}
    \hfill
    % Segunda imagem com etiqueta (b)
    \begin{minipage}{0.48\textwidth}
        \centering
        \includegraphics[width=\textwidth]{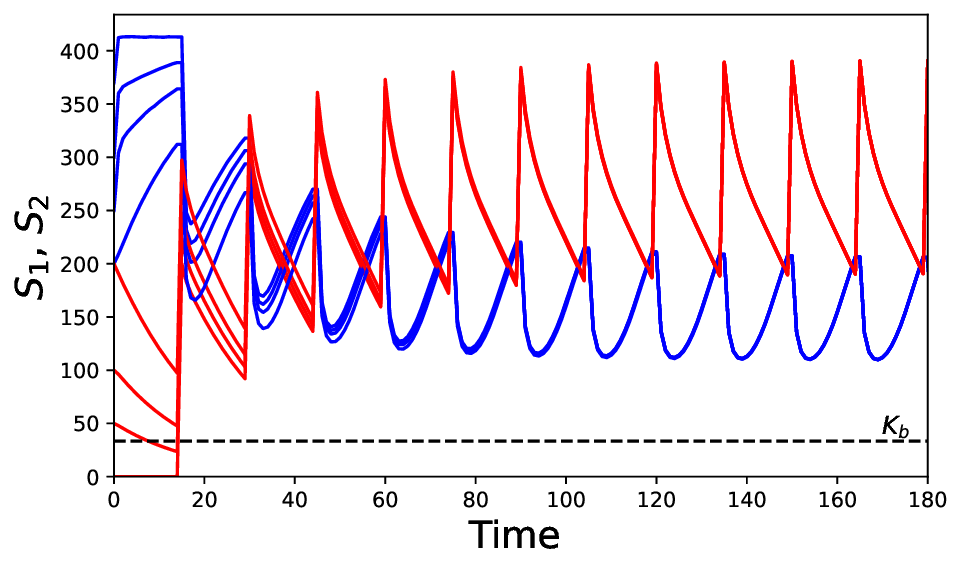}
        \caption*{(b)}
    \end{minipage}
    % Legenda geral para as duas imagens
   \caption{Impulsive solutions of the model~\eqref{eq:equation_1}-\eqref{eq:equation_2} for different initial conditions. In Figure~(a), the release amount is $u_k = 100$ with $\tau = 7$ $(u_k < \eta(7))$, and in Figure~(b), the release amount is $u_k = 200$ with $\tau = 14$ $(u_k < \max\limits_{\tau \geq 0} \eta(\tau))$. In both cases, the $S_1$-free solution does not achieve global stability. The dynamics of wild \textit{Aedes aegypti} female mosquitoes $S_1$ are represented in blue, while those infected with Wolbachia $S_2$ are represented in red.
}\label{fig:1}
\end{figure}

\begin{figure}[!ht]
    \centering
    % Primeira imagem com etiqueta (a)
    \begin{minipage}{0.48\textwidth}
        \centering
        \includegraphics[width=\textwidth]{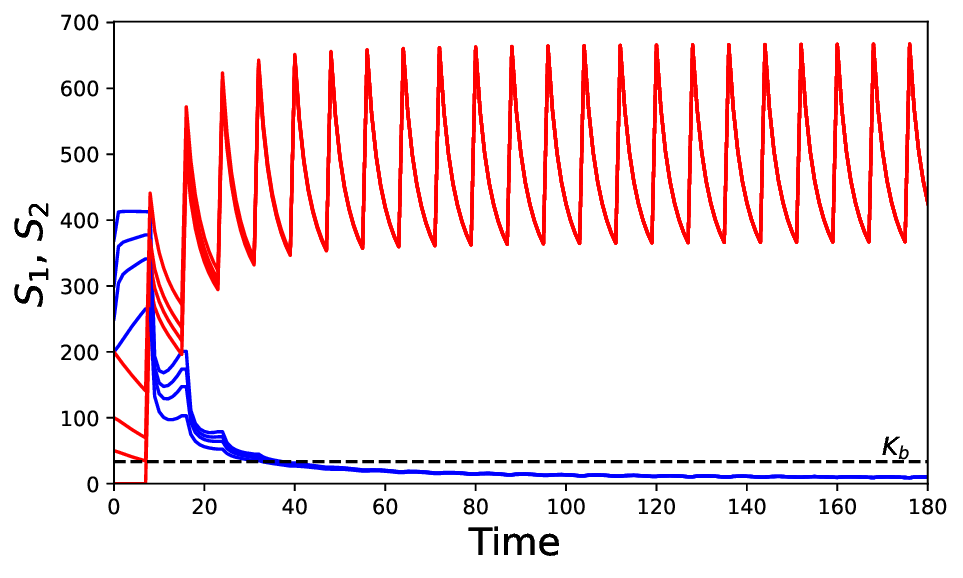}
        \caption*{(c)}
    \end{minipage}
    \hfill
    % Segunda imagem com etiqueta (b)
    \begin{minipage}{0.48\textwidth}
        \centering
        \includegraphics[width=\textwidth]{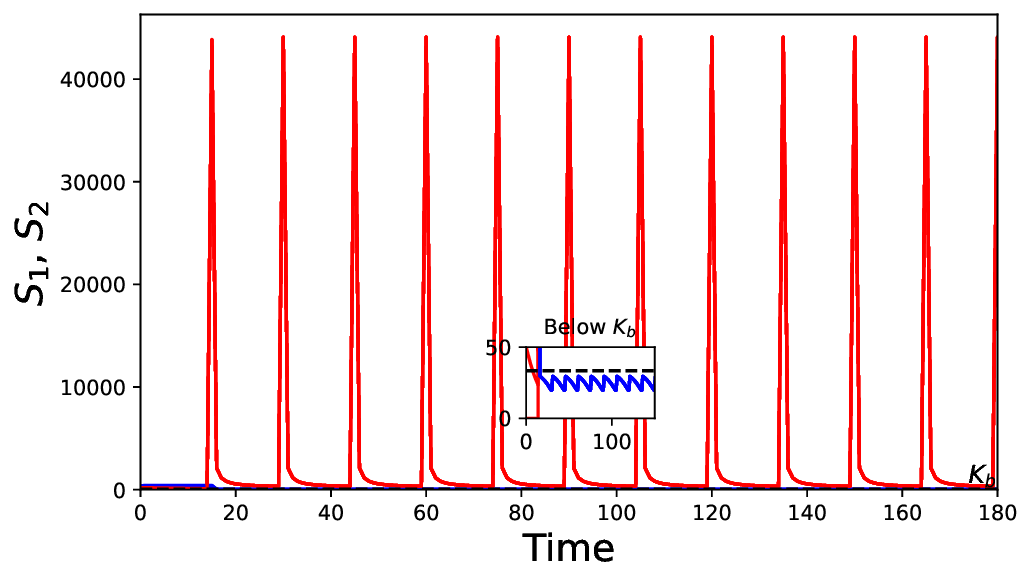}
        \caption*{(d)}
    \end{minipage}
    % Legenda geral para as duas imagens
   \caption{Impulsive solutions of the model~\eqref{eq:equation_1}-\eqref{eq:equation_2} for different initial conditions. In Figure~$(c)$, the release amount is $u_k = 300$ with $\tau = 7$ $(u_k > \eta(7))$, and in Figure~$(d)$, the release amount is $u_k = 43760$ with $\tau = 14$ $(u_k > \max\limits_{\tau \geq 0} \eta(\tau))$. In both cases, the $S_1$-free solution achieves global stability. The dynamics of wild \textit{Aedes aegypti} female mosquitoes $S_1$ are represented in blue, while those infected with Wolbachia $S_2$ are represented in red.}
\label{fig:2}
\end{figure}

\begin{figure}[!ht]
    \centering
    % Primeira imagem com etiqueta (a)
    \begin{minipage}{0.48\textwidth}
        \centering
        \includegraphics[width=\textwidth]{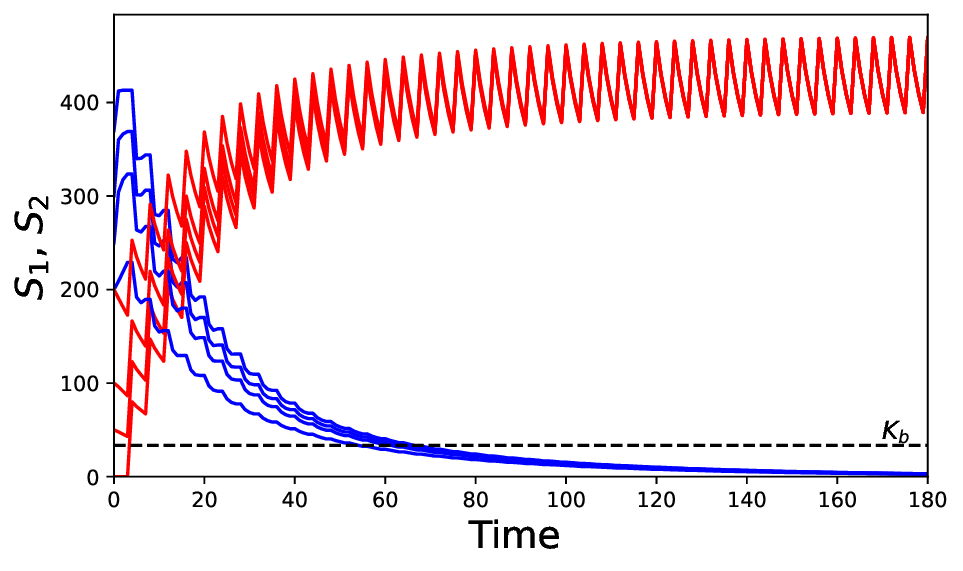}
        \caption*{(e)}
    \end{minipage}
    \hfill
    % Segunda imagem com etiqueta (b)
    \begin{minipage}{0.48\textwidth}
        \centering
        \includegraphics[width=\textwidth]{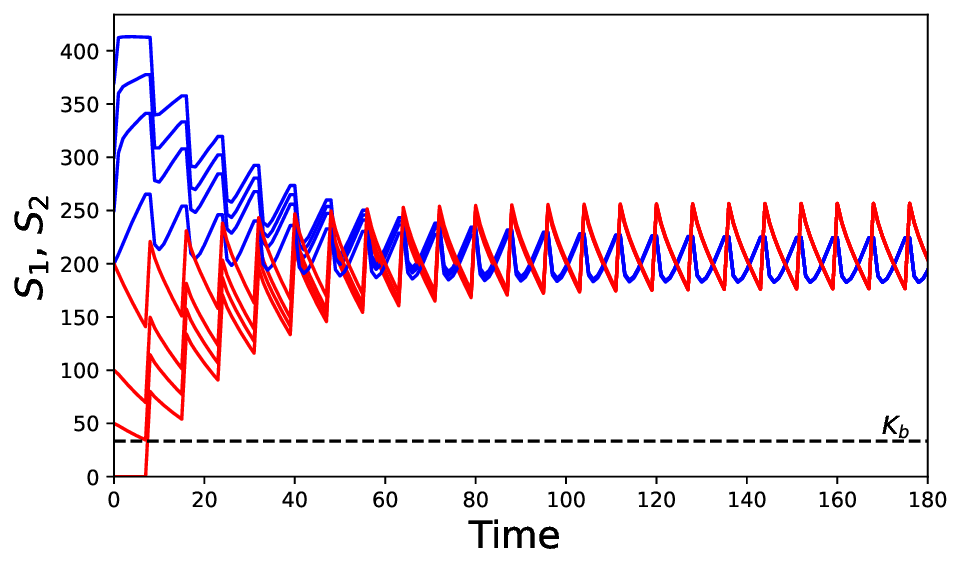}
        \caption*{(f)}
    \end{minipage}
    % Legenda geral para as duas imagens
   \caption{Impulsive solutions of the model \eqref{eq:equation_1}-\eqref{eq:equation_2} for different initial conditions. The release amount is $u_k = 80$, with $\tau = 3$ $(u_k > \eta(3))$ and $\tau = 7$ $(u_k < \eta(7))$, in  Figures~$(e)$ and $(f)$, respectively. In the first case, the $S_1$-free solution achieves global stability, while in the second case it does not. The dynamics of wild \textit{Aedes aegypti} female mosquitoes $S_1$ are shown in blue, and those infected with Wolbachia $S_2$ are shown in red.}

    \label{fig:3}
\end{figure}

\begin{figure}[!ht]
    \centering
    % Primeira imagem com etiqueta (a)
    \begin{minipage}{0.48\textwidth}
        \centering
        \includegraphics[width=\textwidth]{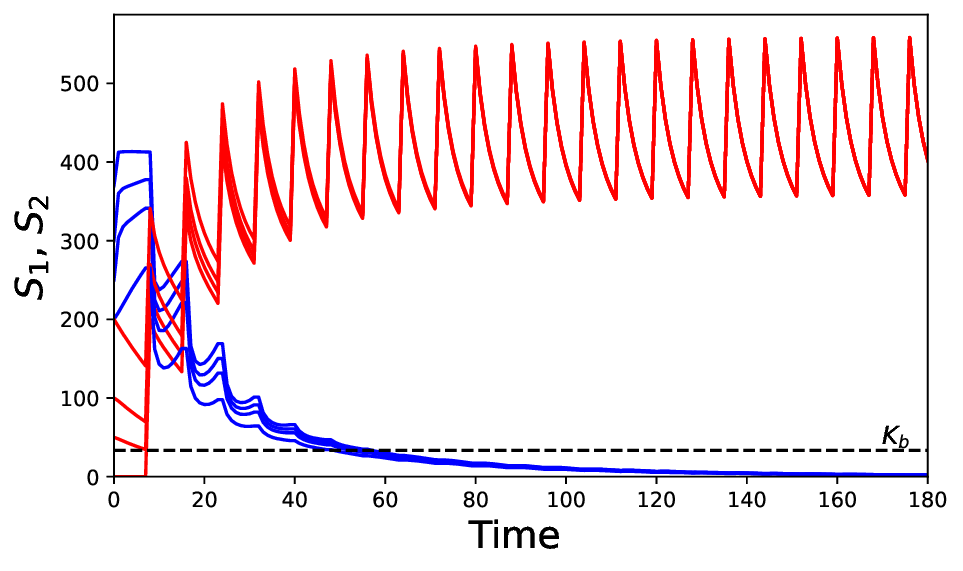}
        \caption*{(g)}
    \end{minipage}
    \hfill
    % Segunda imagem com etiqueta (b)
    \begin{minipage}{0.48\textwidth}
        \centering
        \includegraphics[width=\textwidth]{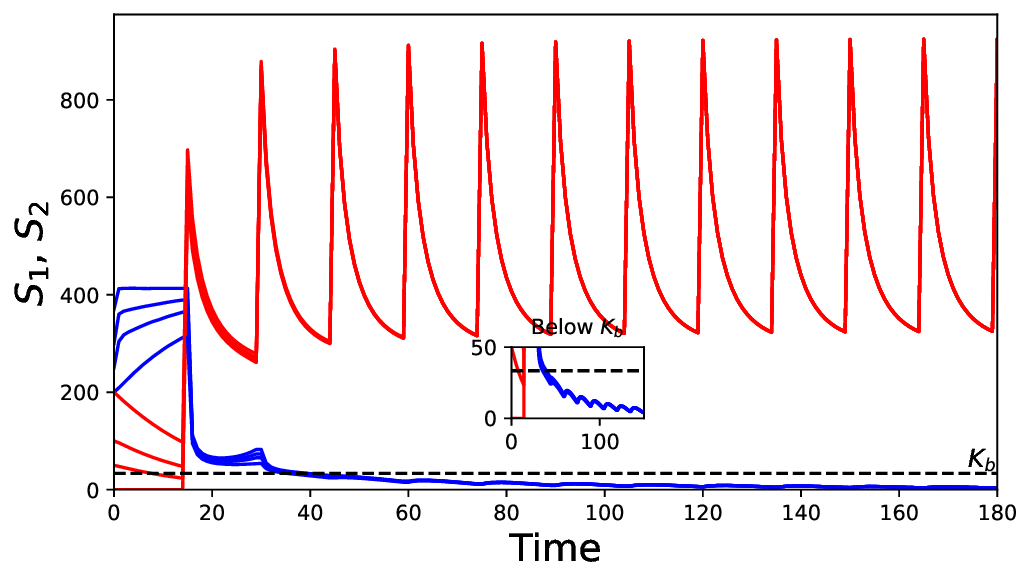}
        \caption*{(h)}
    \end{minipage}
    % Legenda geral para as duas imagens
   \caption{Impulsive solutions of the model \eqref{eq:equation_1}-\eqref{eq:equation_2} for different initial conditions. In Figure~$(g)$, the release amount is $u_k = 200$ with $\tau = 7$ $(u_k < \eta(7))$, and in Figure~$(h)$, the release amount is $u_k = 600$ with $\tau = 14$ $(u_k < \max\limits_{\tau \geq 0} \eta(\tau))$. In both cases, even without satisfying condition~\eqref{eq:equation_14}, the $S_1$-free solution achieves global stability. The dynamics of wild \textit{Aedes aegypti} female mosquitoes $S_1$ are represented in blue, while those infected with Wolbachia $S_2$ are represented in red.
}
    \label{fig:4}
\end{figure}

The simulation results are presented in  \Cref{fig:1,fig:2,fig:3,fig:4}. They illustrate the analysis of the global stability of the wild-female-free solution for different release periods $\tau$.

\Cref{fig:1} (a) and (b) illustrate situations in which global stability is not achieved. This is consistent with the result presented in \Cref{thm:thm_6}, since in both cases the number of infected females $S_2$ released in each period $\tau$ does not satisfy the sufficient condition $u_k > \eta(\tau)$ to ensure global stability. For example, in  \Cref{fig:1} (a), with $\tau = 7$ and a constant release sequence $u_k = 100, \forall k > 0$, the release value does not meet the threshold derived in the previous subsection, as $\eta(7) \approx 300$. A similar situation occurs in \Cref{fig:1} (b), where for $\tau = 14$ and $u_k = 200, \, \forall k > 0$, the release does not satisfy the sufficient condition for the global stability of $(0, \overline{S}_2(t))$ given in equation \eqref{eq:equation_14}.

\Cref{fig:2} (c) and (d) illustrate situations in which global stability is achieved. This occurs when the number of individuals released in each period $\tau$ is increased to satisfy the sufficient condition given by $u_k > \max\limits_{\tau \geq 0} \eta(\tau)$. Specifically, in \Cref{fig:2} (c), for $\tau = 7$, and \Cref{fig:2} (d), for $\tau = 14$, the release is adjusted to constant sequences of $u_k = 300$ and $u_k = 43760, \,\forall k > 0$, respectively, both of which satisfy the stability condition in \eqref{eq:equation_14}.

Observe in the figures how a constant release period $\tau$ and the quantity released $u_k$ directly influence the population dynamics of wild female Aedes aegypti mosquitoes $S_1$ and those infected with Wolbachia $S_2$. For example, in \Cref{fig:3} (e) and \Cref{fig:3} (f), we observe the effect of applying a constant release sequence $u_k = 80,\, \forall k > 0$ with two different release periods: $\tau = 3$ in \Cref{fig:3} (e) and $\tau = 7$ in \Cref{fig:3} (f). 
In \Cref{fig:3} (e), global stability of the solution $(0, \overline{S}_2(t))$ is achieved, as $u_k = 80 > \eta(3) \approx 60, \, \forall k > 0$. By contrast, in \Cref{fig:3} (f), where $u_k = 80 < \eta(7) \approx 300, \, \forall k > 0$, global stability of the solution is not attained.

Based on the simulations performed, the dependence between $\tau$ and $u_k$ underscores the importance of carefully balancing the frequency and quantity of releases, as both parameters directly impact the success in achieving global system stability. Notably, as the interval between releases increases, the required number of individuals to be released also rises. Furthermore, considering the model conditions and the parameters presented in \Cref{tab:tab_1}, we observe that in all implemented scenarios, particularly for $\tau = 14$ with $\eta(14) < 0$, as shown in \Cref{fig:2} (d), the release value $u_k = 43760$ is larger than $\max\limits_{\tau \geq 0} \eta(\tau) \approx 43759.89$, thereby satisfying the condition for global stability. However, \Cref{fig:2} (d) suggests that this amount might exceed what is necessary to eliminate the $S_1$ population.

These observations indicate that the sufficiency of condition \eqref{eq:equation_14} ensures the global stability of the $S_1$-free solution; however, it may be unnecessarily strict. For example, in \Cref{fig:4} (g) and (h), we observe that the global stability of the wild-female-free solution $S_1$ is attained with constant sequences of releases of $u_k = 200$ and $u_k = 600, \, \forall k > 0$, for $\tau = 7$ and $\tau = 14$, respectively.

These simulations provide a deeper understanding of the conditions that either ensure or compromise the system's stability. Furthermore, we observed that it is possible to release fewer infected females than established by the sufficient condition and still achieve global stability. This is because the condition~\eqref{eq:equation_14} is only sufficient to ensure the global stability of the solution $(0, \overline{S}_2(t))$.

Finally, comparing the simulation results with the theoretical analysis, they confirm the theoretical findings for the impulsive model studied. Additionally, we emphasize the need for an optimization strategy for the release process that minimizes intervention costs while ensuring the fixation of Wolbachia-infected females in the target area, a topic that will be further explored in the next subsection.

\subsection{Impulsive optimal control solutions}

After ensuring the existence of a solution for the problem~\eqref{eq:equation_15}-\eqref{eq:equation_18}, we employed \texttt{GEKKO}, an optimization package in Python ~\cite{gekko}. Using the parameters listed in ~\Cref{tab:tab_1} and the constant $C = 1/200$, we simulated four scenarios by varying the release interval $\tau$. As in the simulations from the previous subsection, $S_1$ represents wild \textit{Aedes aegypti} females, and $S_2$ represents females infected with \textit{Wolbachia} bacterium. The values of $\tau$ for each scenario are: $\tau = 7$ (Case~1), $\tau = 14$ (Case~2), $\tau = 21$ (Case~3), and $\tau = 30$ (Case~4).  

For each case, the final intervention time $T$ also varies, taking the values $T = 300$, $T = 180$, $T = 100$, and $T = 70$, respectively. Additionally, for each $T$ and $\tau$, we present the corresponding release amounts $u_k$ over time.

We conducted simulations to analyze the system dynamics and identify optimal release strategies. Comparing different combinations of release intervals $\tau$ and intervention durations $T$ allowed us to balance cost minimization with ensuring the dominance of the infected population $S_2$. This approach offers practical guidance for efficiently achieving biological control objectives.
\begin{table}[h!]
\centering
\begin{minipage}{0.48\textwidth}
\centering
\caption{Numerical optimization results for Case 1 ($\tau = 7$).}
\begin{tabular}{ccc}
\hline
\textbf{$T$} & \textbf{$\sum u_k$} & \textbf{$\min J(u)$} \\ \hline
300                       & 2949.52                                     & 14.75                                           \\ 
180                       & 1243.07                                     & 6.22                                            \\ 
100                       & 1126.62                                     & 5.63                                            \\ 
70                        & 908.66                                      & 4.54                                            \\ \hline
\end{tabular}
\label{tab:tab_2}
\end{minipage}%
\hfill
\begin{minipage}{0.48\textwidth}
\centering
\caption{Numerical optimization results for Case 2 ($\tau = 14$).}
\begin{tabular}{ccc}
\hline
\textbf{$T$} & \textbf{$\sum u_k$} & \textbf{$\min J(u)$} \\ \hline
300                       & 3604.23                                     & 18.02                                           \\ 
180                       & 1677.15                                     & 8.39                                            \\ 
100                       & 1109.34                                     & 5.55                                            \\ 
70                        & 1110.13                                     & 5.55                                            \\ \hline
\end{tabular}
\label{tab:tab_3}
\end{minipage}
\end{table}

\begin{figure}[!ht]
 \begin{subfigure}[t]{1.1\textwidth}
        \centering
        %\textbf{Optimal Control $u^*$ for Case 1}
        \includegraphics[width=\linewidth]{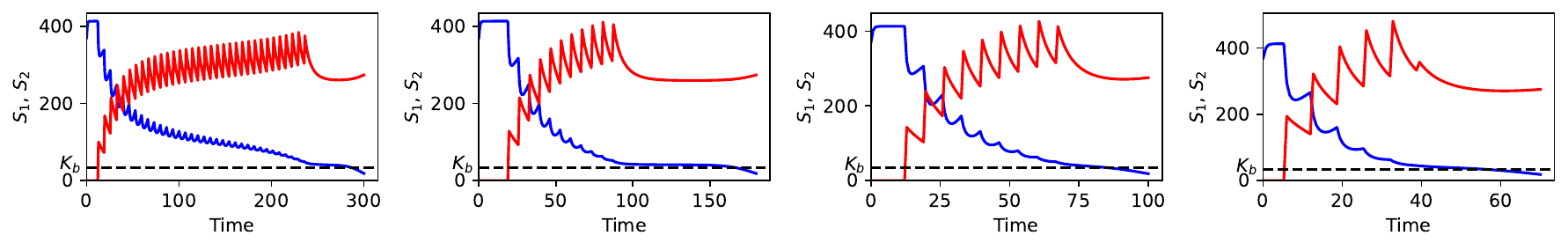}
    \end{subfigure}
    \hfill
    \begin{subfigure}[t]{1.1\textwidth}
        %\textbf{Optimal Trajectory for Case 1}
        \includegraphics[width=\linewidth]{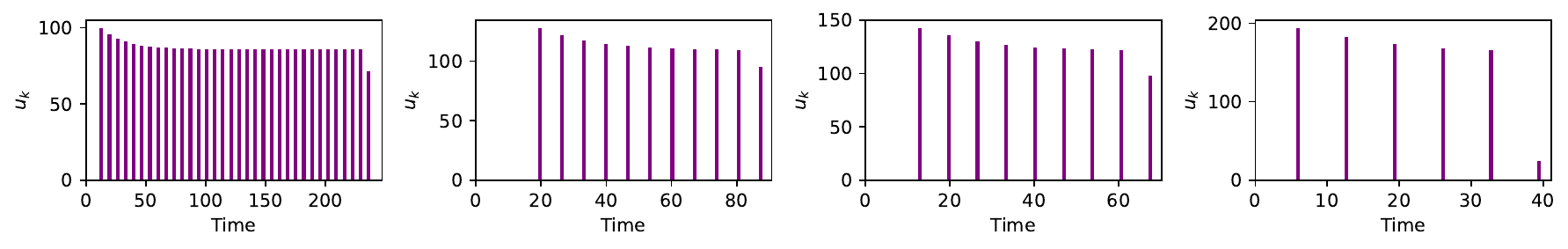}
    \end{subfigure}
    \caption{Optimal trajectories and corresponding controls for Case 1 with initial condition $(K_1, 0)$. $S_1$ in blue; $S_2$ in red. Final times $T$: $300$, $180$, $100$, and $70$ (left to right).}
    \label{fig:5}
\end{figure}

$\bullet$ Case 1: ($\tau = 7$).
In this case, the release interval is set to $\tau = 7$. \Cref{tab:tab_2} summarizes the numerical results, including the total sum of $u_k$ and the minimum value of the functional for each of the four final times $T$ analyzed. \Cref{fig:5} depicts the trajectories of the wild female population $S_1$ and the infected population $S_2$, together with the respective release values $u_k$. For longer intervention horizons, the maximum number of infected females released at a single time step tends to be smaller, as the system has more time to converge to the target state.

\begin{figure}[!ht]
    \begin{subfigure}[t]{1.1\textwidth}
        \centering
        %\textbf{Optimal Control $u^*$ for Case 1}
        \includegraphics[width=\linewidth]{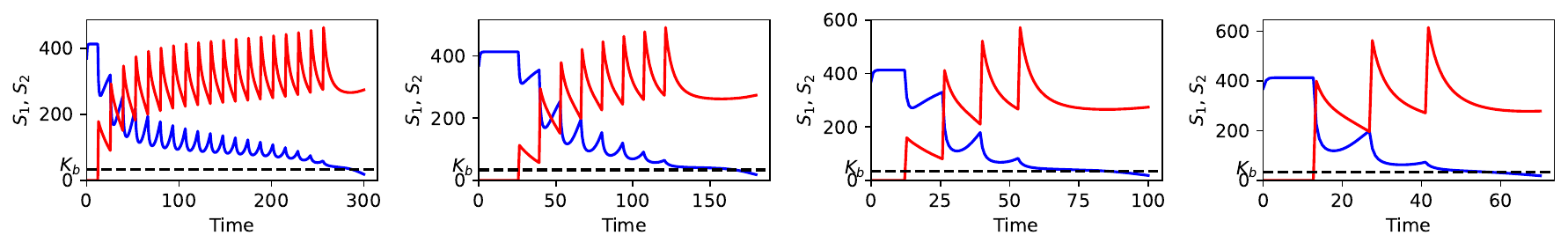}
    \end{subfigure}
    \hfill
    \begin{subfigure}[t]{1.1\textwidth}
        %\textbf{Optimal Trajectory for Case 1}
        \includegraphics[width=\linewidth]{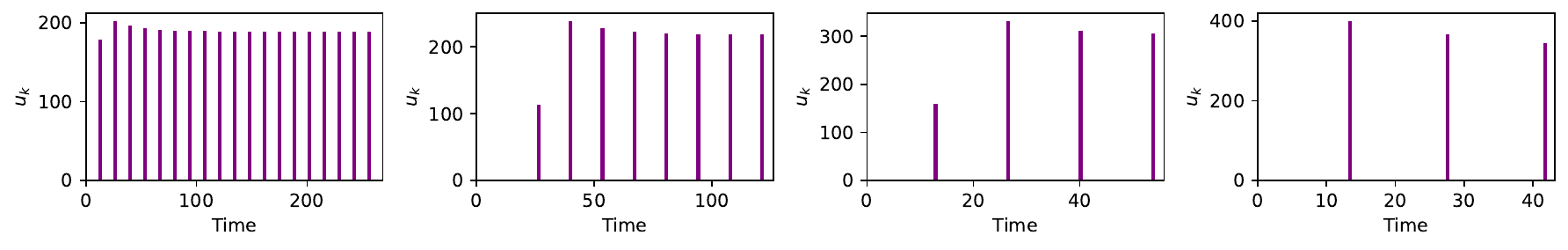}
    \end{subfigure}
    \caption{Optimal trajectories and corresponding controls for Case 2 with initial condition $(K_1, 0)$. $S_1$ in blue; $S_2$ in red. Final times $T$: $300$, $180$, $100$, and $70$ (left to right).}
    \label{fig:6}
\end{figure}

$\bullet$ Case 2: ($\tau = 14$).
Here, the release interval is increased to $\tau = 14$. The corresponding results are presented in \Cref{tab:tab_3} and illustrated in \Cref{fig:6}. In most scenarios, a longer interval between releases required a larger maximum number of infected females per intervention to maintain the wild population below the survival threshold $K_b$. An exception was observed when $T = 100$: despite a higher maximum release per period compared to Case 1, the total number of released females is smaller.

\begin{figure}[H]
   \begin{subfigure}[t]{1.1\textwidth}
        \includegraphics[width=\linewidth]{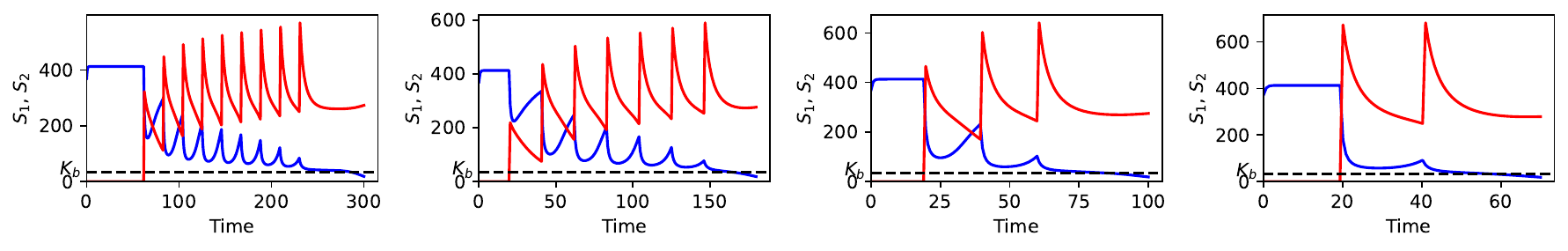}
    \end{subfigure}
    \hfill
    \begin{subfigure}[t]{1.1\textwidth}
        \includegraphics[width=\linewidth]{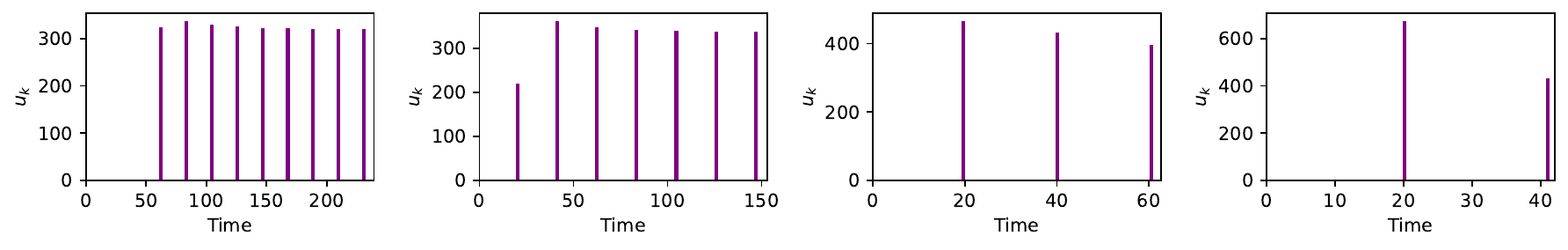}
    \end{subfigure}
    \caption{Optimal trajectories and corresponding controls for Case 3 with initial condition $(K_1, 0)$. $S_1$ in blue; $S_2$ in red. Final times $T$: $300$, $180$, $100$, and $70$ (left to right).}
    \label{fig:7}
 \end{figure}

$\bullet$ Case 3: ($\tau = 21$).
With $\tau = 21$, the system behavior is shown in \Cref{fig:7}, and the numerical outcomes are reported in \Cref{tab:tab_4}. A further increase in the maximum release per period is observed relative to previous cases. The variation in the minimum value of the functional is not monotonic across $\tau$: for $T = 180$ and $T = 100$, the values (11.42 and 6.47, respectively) exceeded those in Cases 1 and 2, whereas for the remaining final times, the trend is reversed.
\begin{table}[h!]
\centering
\begin{minipage}{0.48\textwidth}
\centering
\caption{Numerical optimization results for Case 3 ($\tau = 21$).}
\begin{tabular}{ccc}
\hline
\textbf{$T$} & \textbf{$\sum u_k$} & \textbf{$\min J(u)$} \\ \hline
300                       & 2918.42                                     & 14.59                                           \\ 
180                       & 2283.66                                     & 11.42                                           \\ 
100                       & 1294.69                                     & 6.47                                            \\ 
70                        & 1106.70                                     & 5.53                                            \\ \hline
\end{tabular}
\label{tab:tab_4}
\end{minipage}%
\hfill
\begin{minipage}{0.48\textwidth}
\centering
\caption{Numerical optimization results for Case 4 ($\tau = 30$).}
\begin{tabular}{ccc}
\hline
\textbf{$T$} & \textbf{$\sum u_k$} & \textbf{$\min J(u)$} \\ \hline
300                       & 3788.17                                     & 18.94                                           \\ 
180                       & 2081.34                                     & 10.41                                           \\ 
100                       & 1256.42                                     & 6.28                                            \\ 
70                        & 873.47                                      & 4.37                                            \\ \hline
\end{tabular}
\label{tab:tab_5}
\end{minipage}
\end{table}
    
\begin{figure}[H]
    \begin{subfigure}[t]{1.1\textwidth}
        \includegraphics[width=\linewidth]{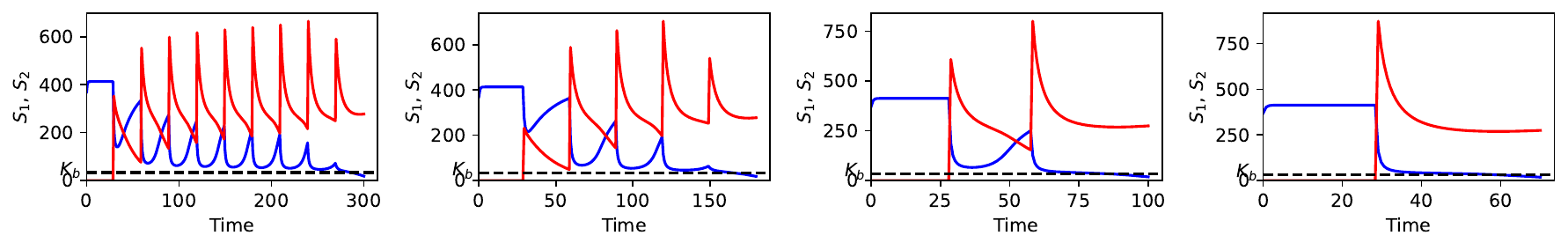}
    \end{subfigure}
    \hfill
    \begin{subfigure}[t]{1.1\textwidth}
        \includegraphics[width=\linewidth]{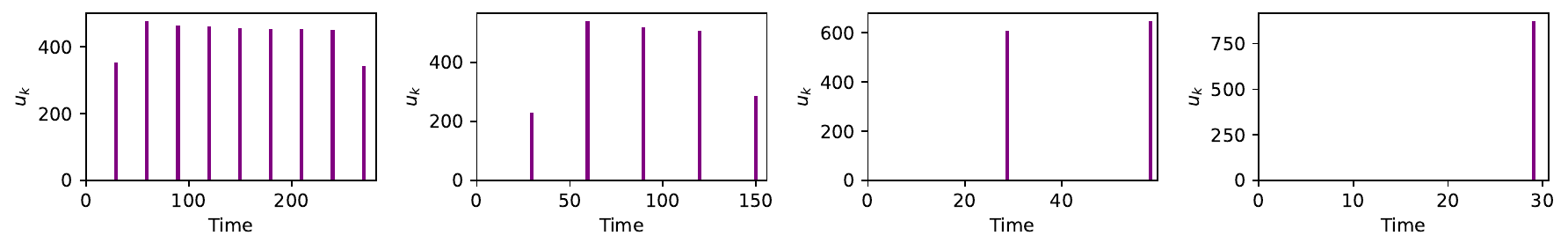}
    \end{subfigure}
    \caption{Optimal trajectories and corresponding controls for Case 4 with initial condition $(K_1, 0)$. $S_1$ in blue; $S_2$ in red. Final times $T$: $300$, $180$, $100$, and $70$ (left to right).}
    \label{fig:8}
    \end{figure}

$\bullet$ Case 4: ($\tau = 30$).
In this case, the release interval was set to $\tau = 30$. The results are given in \Cref{tab:tab_5}, and the corresponding trajectories are depicted in \Cref{fig:8}. For $T = 300$, this configuration requires the most significant overall number of infected females, while Case 3 yields the smallest total release. Notably, for $T = 70$, a single release of $u_k = 873.47$ was sufficient to bring the wild female population to its survival threshold, ensuring the fixation of $S_2$ in the target area.
It is important to note that these results rely on the assumption that all model parameters remain constant over time. Parameters such as birth and death rates vary due to environmental factors. For mosquitoes, for instance, these rates depend on humidity, temperature, and other seasonal conditions. Thus, although the model predicts that a single release (Case 4, $T = 70$) would suffice under the assumed conditions, it remains a simplified mathematical representation. This reinforces the need for multiple simulations with different release strategies.

The simulations reveal a strong relationship between the impulse period $\tau$, the release amount $u_k$, and the chosen final time $T$. Increasing $\tau$ raises the maximum release per intervention, but the total number of infected females released does not follow a monotonic trend. The interaction between $\tau$ and $T$ determines efficiency: for $T = 100$, $\tau = 14$ required fewer total releases than $\tau = 7$, despite a higher per-event peak; for $T = 180$, the highest total release occurred at $\tau = 21$, and the lowest at $\tau = 7$. This indicates that neither parameter alone determines performance; rather, their combination defines the optimal strategy.

These findings highlight the value of exploring multiple $(\tau, T)$ configurations to identify strategies best suited to specific constraints, such as available time or the number of infected females for release. Overall, the simulations confirm the effectiveness of the optimal control approach for Wolbachia-infected female releases.

%\section{Discussions}

\section{Final Remarks}
\label{sec:conclusions}
This study presents a mathematical framework to investigate the dynamics of two competing species under impulsive release strategies. By defining an impulsive differential equation model, we demonstrate the potential of periodic interventions in ecological systems, particularly in reducing the population of one species $S_1$ through strategic competition with another species $S_2$.

Theoretical analysis confirmed that the model is well-posed, ensuring the existence and uniqueness of an $S_1$-free solution under specific conditions. Moreover, we derived a sufficient condition for impulsive releases of $S_2$ that guarantees the global stability of the desired equilibrium, providing a solid foundation for ecological control strategies.

In this work, we adopt impulsive releases as an alternative to the continuous controls often considered in the literature, better reflecting realistic intervention settings. The impulsive nature of the releases justifies the use of numerical methods to determine optimal release magnitudes, ensuring that the resulting strategy is both effective and feasible.

This work highlights the relevance of mathematical modeling and control theory in tackling ecological challenges. Future research would extend these methods to more complex environmental systems, considering time-varying parameters, assessing how environmental factors influence control effectiveness, and analyzing how parameter variability affects stability and optimality. Additional directions include evaluating the interplay between different control methods and exploring model adaptations that better capture realistic ecological dynamics. Such developments could significantly enhance the applicability and robustness of impulsive intervention strategies in ecological control. Moreover, extending the problem to a framework involving fractional derivatives represents an interesting avenue that could be explored in future studies.

\appendix
\section{Proof of Theorem \ref{thm:thm_3}} 
\label{appendix:A}
\begin{proof}
First, we will show that $\overline{Z}_2(t)$ is the unique positive $\tau$-periodic solution of system \eqref{eq:equation_6}. Let
\begin{equation*}
    Z_2(t) = \dfrac{c K_2 e^{r_2 t}}{c e^{r_2 t} - 1},
\end{equation*}
be the solution of
\begin{equation*}
    \dfrac{dZ_2}{dt} = Z_2\left(\psi_2 - \dfrac{r_2}{K_2}Z_2\right) - \delta_2 Z_2, \quad t \neq k\tau, \quad k \geq 0,
\end{equation*}
which corresponds to the first equation of system \eqref{eq:equation_6}, where $c \in \mathbb{R}$ is a constant to be determined and is associated with the initial conditions of the problem.

For $t = k\tau$, $k \geq 0$, let $Z_2(k\tau^+)$ denote the initial value at time $k\tau$. Then,
\begin{equation*}
    c = \dfrac{Z_2(k\tau^+) e^{-r_2 k\tau}}{Z_2(k\tau^+) - K_2},
\end{equation*}
and the solution becomes
\begin{equation}\label{eq:equation_19}
    Z_2(t) = \dfrac{K_2 Z_2(k\tau^+) e^{r_2 (t-k\tau)}}{Z_2(k\tau^+)\left(e^{r_2 (t-k\tau)} - 1\right) + K_2}, \quad k\tau < t \leq (k+1)\tau, \quad k \geq 0,
\end{equation}
which represents the solution of system \eqref{eq:equation_6} between the pulses.

At the impulsive moments, when $t = (k+1)\tau$, $k \geq 0$, from the second equation of \eqref{eq:equation_6}, with $u_k \in U$, we obtain the following difference equation:
\begin{eqnarray}\label{eq:equation_20}\nonumber
    Z_2((k+1)\tau^+) &= Z_2((k+1)\tau^-) + u_k \\
    &= \dfrac{K_2 Z_2(k\tau^+) e^{r_2 \tau}}{Z_2(k\tau^+)\left(e^{r_2 \tau} - 1\right) + K_2} + u_k,
\end{eqnarray}
which defines a recursive relation between $Z_2((k+1)\tau^+)$ and $Z_2(k\tau^+)$. We can rewrite it as
\begin{equation}\label{eq:equation_21}
    Z_2^{k+1} = \dfrac{K_2 Z_2^k e^{r_2 \tau}}{Z_2^k\left(e^{r_2 \tau} - 1\right) + K_2} + u_k,
\end{equation}
where $Z_2^k = Z_2(k\tau^+)$.

Define
\begin{equation*}
    h(Z_2) = \dfrac{K_2 Z_2 e^{r_2 \tau}}{Z_2\left(e^{r_2 \tau} - 1\right) + K_2} + u_k,
\end{equation*}
and observe that \eqref{eq:equation_21} has a single positive equilibrium point. Indeed, setting
\begin{equation*}
    h(Z_2) = Z_2,
\end{equation*}
gives the condition that determines the equilibrium point of the recursive equation \eqref{eq:equation_21}, which is $\tau$-periodic. Performing algebraic manipulations, we obtain the quadratic equation
\begin{equation*}
    Z_2^2 - (u_k + K_2)Z_2 - \dfrac{u_k K_2}{e^{r_2 \tau} - 1} = 0,
\end{equation*}
which has two real roots, but only one of them is positive. It is given by
\begin{equation*}
    Z_2^+ = \dfrac{1}{2} \left[(u_k + K_2) + \sqrt{(u_k + K_2)^2 + 4 \dfrac{u_k K_2}{e^{r_2 \tau} - 1}} \right], \quad k \geq 0.
\end{equation*}
Then, substituting $Z_2^+$ into \eqref{eq:equation_19}, we obtain the corresponding $\tau$-periodic solution:
\begin{equation*}
    \overline{Z}_2(t) = \dfrac{K_2 Z_2^+ e^{r_2 (t-k\tau)}}{Z_2^+\left(e^{r_2 (t-k\tau)} - 1\right) + K_2}, \quad k\tau < t \leq (k+1)\tau, \quad k \geq 0,
\end{equation*}
which is the unique positive $\tau$-periodic solution of system \eqref{eq:equation_6}.

Now we show that $\overline{Z}_2(t)$ is globally asymptotically stable. Since $h(Z_2)$ is continuous and differentiable for $Z_2 > 0$, we compute its derivative:
$$
h'(Z_2) = \dfrac{\partial}{\partial Z_2} \left( \dfrac{K_2 Z_2 e^{r_2 \tau}}{Z_2(e^{r_2 \tau} - 1) + K_2} + u_k \right)
= \dfrac{K_2^2 e^{r_2 \tau}(e^{r_2 \tau} - 1)}{\left(Z_2(e^{r_2 \tau} - 1) + K_2\right)^2} > 0,
$$
for all $Z_2 > 0$. Hence, $h$ is strictly increasing.

We have already shown that $Z_2^+$ is the unique positive fixed point of $h$, that is,
$$
h(Z_2^+) = Z_2^+,
$$
and there is no other solution to $h(Z_2) = Z_2$ with $Z_2 > 0$.

Since $h$ is strictly increasing and $Z_2^+$ is its unique positive fixed point, it follows that:
\begin{itemize}
    \item If $0 < Z_2 < Z_2^+$, then $h(Z_2) < h(Z_2^+) = Z_2^+$, so $Z_2 < h(Z_2) < Z_2^+$.
    \item If $Z_2 > Z_2^+$, then $h(Z_2) > h(Z_2^+) = Z_2^+$, so $Z_2 > h(Z_2) > Z_2^+$.
\end{itemize}

These inequalities show that the sequence $\{Z_2^k\}$ defined by the recursive equation \eqref{eq:equation_20} is monotonic and converges to $Z_2^+$ for any initial condition $Z_2^0 > 0$. Therefore, $Z_2^+$ is a globally asymptotically stable fixed point of the discrete dynamical system.

As a consequence, the corresponding $\tau$-periodic solution $\overline{Z}_2(t)$ of the impulsive system \eqref{eq:equation_6} is globally asymptotically stable.

\end{proof}

\section{Proof of Lemma \ref{lem:lem_1}} 
\label{appendix:B}
\begin{proof}
From the first equation of the system \eqref{eq:equation_1}–\eqref{eq:equation_2}, we have
\begin{equation*}
    \dfrac{dS_1}{dt}(t) = S_1(t)\left(\psi_1-\dfrac{r_1}{K_1}(S_1(t)+S_2(t))\right) \left(\dfrac{S_1(t)}{K_0} - 1 \right) - \delta_1 S_1(t),
\end{equation*}
which does not involve impulsive releases. The presence of the critical depensation term $\left(\frac{S_1(t)}{K_0} - 1\right)$, associated with the concept of minimum viable population size (MVPS) \cite{Campo2017}, requires splitting the analysis into two cases. Here, $K_0$ plays the role of a critical threshold: if the population falls below this value, it tends to extinction; otherwise, it may persist. Analyzing both cases allows us to conclude the boundedness of $S_1(t)$.

\begin{itemize}
    \item[\textbf{a.}] If \(S_1(t) \leq K_0\), then \(\left(\frac{S_1(t)}{K_0} - 1 \right) \leq 0\). Thus, \(K_0\) is an upper bound for \(S_1(t)\) for all \(t \geq 0\).

    \item[\textbf{b.}] If \(S_1(t) > K_0\), then \(\left(\frac{S_1(t)}{K_0} - 1 \right) > 0\). We then have
\begin{align*}
   \dfrac{dS_1}{dt}(t) &= S_1(t)\left[\left(\psi_1 -\dfrac{r_1}{K_1}S_1(t)\right)\left(\dfrac{S_1(t)}{K_0} - 1 \right)-\dfrac{r_1}{K_1}S_2(t)\left(\dfrac{S_1(t)}{K_0} - 1 \right)-\delta_1\right] \nonumber\\
   &\leq S_1(t)\left[\left(\psi_1 -\dfrac{r_1}{K_1}S_1(t)\right)\left(\dfrac{S_1(t)}{K_0} - 1 \right)-\delta_1\right],
\end{align*}
and we can consider the comparison differential equation
\begin{align}\label{eq:equation_22}
\begin{matrix}
\dfrac{dy}{dt}(t) = & y(t)\left[\left(\psi_1 -\dfrac{r_1}{K_1}y(t)\right)\left(\dfrac{y(t)}{K_0} - 1 \right)-\delta_1\right], \\
y(0) = & \hspace{-5.5cm} S_1(0),
\end{matrix}
\end{align}
which has three equilibrium points: \(y = 0\), and the two positive solutions of
\begin{equation}\label{eq:equation_23}
    \left(\psi_1 -\dfrac{r_1}{K_1}y\right)\left(\dfrac{y}{K_0} - 1 \right)-\delta_1 = 0.
\end{equation}
Expanding \eqref{eq:equation_23}, we obtain the quadratic equation
\begin{equation*}
-\dfrac{r_1}{K_1K_0}y^2 + \dfrac{\psi_1K_1 + r_1K_0}{K_1K_0}y - (\delta_1 + \psi_1) = 0.
\end{equation*}
Solving this, we find
\begin{equation*}
\Delta = (\psi_1K_1 + r_1K_0)^2 - 4r_1K_0K_1(\delta_1 + \psi_1),
\end{equation*}
and the roots
\begin{equation*}
     y_1 = \dfrac{\psi_1K_1 + r_1K_0 - \sqrt{\Delta}}{2r_1}, \quad 
     y_2 = \dfrac{\psi_1K_1 + r_1K_0 + \sqrt{\Delta}}{2r_1}.
\end{equation*}
Following the notation in \cite{Campo2017}, we define \(y_1 = K_b\) and \(y_2 = K_*\), with \(0 < K_b < K_*\). We now analyze the sign of the right-hand side of \eqref{eq:equation_22} in the intervals between the equilibria:
\begin{equation*}
    \left \{
    \begin{aligned}
        \dfrac{dy}{dt} < 0, &\quad \text{if } 0 < y < K_b, \\
        \dfrac{dy}{dt} > 0, &\quad \text{if } K_b < y < K_*, \\
        \dfrac{dy}{dt} < 0, &\quad \text{if } y > K_*.
    \end{aligned}
    \right.
\end{equation*}
Therefore, for the initial condition $y(0) = S_1(0)$, we have:
\begin{itemize}
    \item[i)] If $0 \leq S_1(0) < K_b$, then $y(t)$ decreases to $0$ in finite time as $t$ increases;
    \item[ii)] If $K_b < S_1(0) < K_*$, then $y(t)$ increases toward $K_*$ as $t \to \infty$;
    \item[iii)] If $S_1(0) > K_*$, then $y(t)$ decreases toward $K_*$ as $t \to \infty$.
\end{itemize}
As a result, following \cite{Campo2017}, given an initial population of species $S_1$ large enough to survive (i.e., $S_1(0) > K_b$), the threshold $K_b$ separates extinction from persistence, and $K_*$ represents the stable steady-state population. Since the comparison equation satisfies $\frac{dS_1}{dt}(t) \leq \frac{dy}{dt}(t)$ and $S_1(0) = y(0)$, by the Comparison Theorem we conclude that

$$S_1(t) \leq y(t) \leq K_*, \quad \text{for all } t \geq 0.$$
\end{itemize}
Combining both cases (a) and (b), we conclude that $S_1(t)$ is bounded for all $t \geq 0$. In particular,

$$S_1(t) \leq M_1 := \max\{K_*, S_1(0)\}, \quad \text{for all } t \geq 0.$$

\end{proof}

\section{Proof of Theorem \ref{prop:prop_3}} 
\label{appendix:C}
\begin{proof}
By hypothesis, the initial conditions are non-negative. By Proposition \ref{prop:prop_1}, we have that $S_1(t)$ and $S_2(t)$ are lower bounded by zero for all $t \geq 0$, and by Lemma \ref{lem:lem_1}, $M_1 > 0$ is an upper bound for $S_1(t)$. We will now show that $S_2(t)$ is also upper-bounded.
    To do this, consider the second and fourth equations of system \eqref{eq:equation_1}--\eqref{eq:equation_2}, from which we obtain:
    \begin{equation*}
    \dfrac{dS_2}{dt}(t) \leq \dfrac{dZ_2}{dt}(t), \quad \text{with } S_2(0) = Z_2(0),
    \end{equation*}
    where $Z_2$ satisfies \eqref{eq:equation_6}. The function $Z_2(t)$ is bounded, since for $t \neq k\tau$, the solution of the continuous ODE is bounded by $\max\{K_2, Z_2(0)\}$, where $K_2$ represents the environmental carrying capacity of $S_2$. At the impulsive instants $t = k\tau$, a term $u_k \in U$ is added, which is upper bounded by $u_{\max}$. Thus, immediately after the impulse, we have $Z_2(t^+) \leq \max\{K_2, S_2(0)\} + u_{\max}$. Therefore, $Z_2(t)$ remains bounded for all time.
    By Theorem \ref{thm:thm_2}, we obtain:
    \begin{equation*}
    S_2(t) \leq Z_2(t).
    \end{equation*}
    Let $M_2 := \max\{K_2, S_2(0)\} + u_{\max}$. Then, $S_2(t) \leq M_2$ for all $t \geq 0$.

To show that the solutions are uniformly bounded, consider the function $V(t) = S_1(t) + S_2(t)$. Then $V(t) \in \mathcal{V}_0$, and for some $\lambda > 0$ and $k\tau \leq t \leq (k+1)\tau$, we have:
\begin{align*}
D^+V(t) + \lambda V(t) &= D^+S_1(t) + D^+S_2(t) + \lambda(S_1(t) + S_2(t))\\\nonumber
&= S_1(t)\left(\psi_1 - \dfrac{r_1}{K_1}(S_1(t) + S_2(t))\right)\left(\dfrac{S_1(t)}{K_0} - 1\right) - \delta_1 S_1(t)\\\nonumber
&\quad + S_2(t)\left(\psi_2 - \dfrac{r_2}{K_2}(S_1(t) + S_2(t))\right) - \delta_2 S_2(t) + \lambda(S_1(t) + S_2(t))\\
&\leq (r_1 + \lambda)M_1 + (r_2 + \lambda)M_2 := M_3.
\end{align*}

At the impulsive times $t = k\tau$, we have $V(k\tau^+) = V(k\tau) + u_k$, for $u_k \in U$. Then, by Lemma 2.2 in \cite{Bainov1993}, we obtain:
\begin{align*}\label{eq:equation_22}\nonumber
V(t) &\leq V(0) e^{-\lambda t} + \int_0^t M_3 e^{-\lambda (t-s)},ds + \sum_{0 \leq k\tau \leq t} u_k e^{-\lambda(t - k\tau)}\\
&\leq V(0) e^{-\lambda t} + \dfrac{M_3}{\lambda}(1 - e^{-\lambda t}) + \sum_{0 \leq k\tau \leq t} u_{\max} e^{-\lambda(t - k\tau)}.
\end{align*}

Thus, as $t \to \infty$, we obtain:
\begin{equation*}
V(t) \leq \dfrac{M_3}{\lambda} + u_{\max} \dfrac{e^{\lambda \tau}}{e^{\lambda \tau} - 1}.
\end{equation*}

In this way, $V(t)$ is uniformly bounded, and by its definition, each positive solution $S_1(t)$ and $S_2(t)$ of system \eqref{eq:equation_1}--\eqref{eq:equation_2} is also uniformly bounded.
\end{proof}

\section*{Acknowledgments}
The authors acknowledge the use of Artificial Intelligence (OpenAI’s ChatGPT) to review the English language, correct orthography, and improve the readability of the manuscript.

J.C.S.A. thanks CAPES (Finance Code 001) for the scholarship. C.E.S. acknowledges support by FEEI-PROCIENCIA-CONACYT-PRONII-STIC and AmSud project BIO-CIVIP-AMSU99.

%Bibliography
\bibliographystyle{unsrt}  
\bibliography{references}

\end{document}